\documentclass{article}
                                                        
\usepackage{fullpage}
\usepackage{amsthm} 
\usepackage{enumitem} 
\usepackage{hyperref}
\usepackage{amssymb} 
\usepackage{amsmath} 
\usepackage{bm} 
\usepackage{changepage}  
\usepackage{mathtools}   
\usepackage{subfig}
\usepackage[final]{changes}

\newcommand{\azul}[1]{{#1}}

\usepackage{cite}

\usepackage{tikz}
\usetikzlibrary{shapes,arrows}
\usetikzlibrary{arrows.meta}
\usepackage{pgfplotstable}  

\usepackage{booktabs}

\newtheorem{prop}{Proposition} 

\theoremstyle{definition} 

\newtheorem{pro}{Problem}
\newtheorem{ass}{Assumption}
\newtheorem{alg}{Algorithm}

\title{\raggedright \textbf{\added{Automatic and }Computationally Efficient Alignment in \added{Fan- and} Cone-beam Tomography}  \\ \vspace{20pt} \normalsize Patricio Guerrero$^{1,*}$, Simon Bellens$^{1,2}$, Ricardo Santander$^{1,2}$, Wim Dewulf$^1$ \\
$^1$Department of Mechanical Engineering, KU Leuven, Celestijnenlaan 300, 3001 Leuven, Belgium \\
$^2$Materialise NV, Technologielaan 15, Leuven, Belgium\\ \medskip
$^*$Corresponding author: Patricio Guerrero, patricio.guerrero@kuleuven.be
\vspace{-30pt}}
\author{}
\date{}

\DeclareMathOperator*{\argmin}{argmin}
\DeclareMathOperator*{\argmax}{argmax}
\DeclareMathOperator*{\shift}{shift}

\newcommand{\R}{\mathbb{R}}
\newcommand{\s}{\mathbb{S}}
\newcommand{\cone}{\mathcal{C}}
\newcommand{\fan}{\mathcal{Q}} 
\newcommand{\g}{\tilde g}
\newcommand{\gb}{\mathbf{g}}
\newcommand{\p}{\mathbf{p}}
\newcommand{\q}{\replaced{\mathbf{w}}{\mathbf{q}}}
\newcommand{\sd}{\mathbf{z}}
\newcommand{\ran}{\mathfrak R}
\newcommand{\h}{h^*}

\newcommand{\et}{\eta^*}
\newcommand{\w}{\Phi} 
\newcommand{\z}{\Psi}  
\newcommand{\x}{\Upsilon} 
\newcommand{\ea}{\emph{et al}}
\newcommand{\eg}{\emph{e.g.}}
\newcommand{\ie}{\emph{i.e.}}
\newcommand{\intpi}{\int\limits_0^{2\pi}}

\begin{document}

\maketitle

\begin{adjustwidth}{30pt}{30pt}
\noindent\textbf{\small Abstract.} 
{\small This work is concerned with fan- and cone-beam computed tomography with circular source trajectory, where the reconstruction inverse problem requires an accurate knowledge of source, detector and rotational axis relative positions and orientations. We address this additional inverse problem as a preceding step of the reconstruction process directly from the acquired projections. In the cone-beam case, we present a method that estimates both the detector shift (orthogonal to both focal and rotational axes) and the in-plane detector rotation (over the focal axis) based on the  variable projection optimization approach. In addition and for the fan-beam case, two new strategies with low computational cost are presented to estimate the detector shift based on a fan-beam symmetry condition. The methods are validated with simulated and \replaced{experimental}{real} industrial tomographic data with code examples available for both fan- and cone-beam geometries. \\
\textbf{keywords}: cone-beam, fan-beam, computed tomography, alignment  \\ \\
\textbf{Published in}: \emph{IEEE Trans. Comput. Imaging} (2024), \href{https://doi.org/10.1109/TCI.2024.3396385}{doi.org/10.1109/TCI.2024.3396385}}
\end{adjustwidth}

\section{Introduction}

\begin{figure}[!t]
\centering
\begin{tikzpicture}[scale=0.9]

  \begin{scope}[shift={(-4,0,0)}]
  \fill[gray!10] (0,2,-2) -- (0,2,2) -- (0,-2,2) -- (0,-2,-2) -- cycle; 
  \end{scope}


  \fill[blue!15] (-4,1.5,0)--(4,0)--(-4,-1.5,0)--cycle;
  
  \begin{scope}[shift={(-4,0,0)}]
  \fill[blue!20] circle (0.6 and 1.5);
  \end{scope}

  \draw (4,-0.1,0) node[anchor=north east] {source};

  \coordinate (O) at (0,0,0);
  \draw[>=latex,->] (0,0,0) -- (-1,0,0) node[anchor=south]{$z$};
  \draw[>=latex,->] (0,0,0) -- (0,1,0) node[anchor=east]{$y$};
  \draw[>=latex,->] (0,0,0) -- (0,0,-1) node[anchor=south]{$x$};

  \begin{scope}[shift={(-4,0,0)}]
  \draw[gray] (0,2,-2) -- (0,2,2) -- (0,-2,2) -- (0,-2,-2) -- cycle; 
  
  \draw[>=latex,->] (0,0,0) -- (0,0,-1) node[anchor=south]{$u$};
  \draw[>=latex,->] (0,0,0) -- (0,1,0) node[anchor=east]{$v$};
  \draw[dotted] (0,2,0) -- (0,-2,0) node[anchor=west]{detector};
  \draw[dotted] (0,0,-2) -- (0,0,2); 
  \draw[<->] (0,2.2,-0.5) -- (0,2.2,0.5); 
  \draw (0,2.2,0.2) node[anchor=south]  {$h$};

  \end{scope}


  \draw[dotted] (4,0,0) -- (-4,0,0); 
  \draw[dotted] (0,2,0) -- (0,-1,0) node[anchor=north]{rotational axis}; 
  \draw[>=latex,->] (-0.3,1.5,0) arc (-150:150:0.3 and 0.1) node[anchor=south]{$\beta$};
  \draw[<->] (-2,0,0) arc (0:240:0.1 and 0.3);
  \draw[] (-2.1,0,0) node[anchor=north west]{$\eta$};

\end{tikzpicture}

\caption{Cone-beam parameterization. The alignment variables to estimate are $\{h,\eta\}$.}
\label{fig.conebeam}
\end{figure}
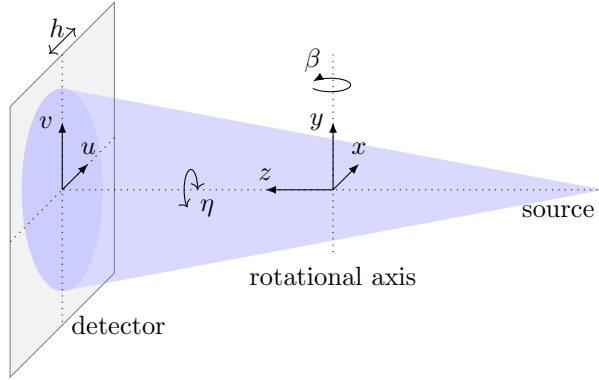

\noindent
This work studies the alignment, or geometry estimation, problem in \azul{fan- and }cone-beam computed tomography (CT) as a preceding step to tomographic reconstruction. The goal in CT is to obtain a three-dimensional (3D) representation, or 3D image, of a target object from a set of 2D measured (X-ray) attenuation images, or projections, under different views of such object. From the acquired projections, a discrete representation of the object is reconstructed typically by means of analytical or iterative numerical strategies \cite{kak}, \azul{\eg\ the industry standard \emph{Feldkamp-Davis-Kress} (FDK) algorithm} \cite{natteredwubbeling}. Cone-beam CT refers here when the projections are taken with X-rays (or other penetrating electromagnetic radiation) that emanate from an external point source forming a cone towards a 2D detector positioned on the other side of the object, as illustrated in Figure \ref{fig.conebeam}. Additionally, we consider a circular trajectory of the source and detector around the object, or equivalently, rotations of the object over a rotational axis parallel to the detector. \azul{This results on fixed relative positions and orientations of source and detector across tomographic rotations}. Such situations are usual in \eg,  industrial or medical CT scanners and are also starting to be implemented in synchrotron radiation facilities.   

Standard reconstruction techniques assume a perfect knowledge of the system geometry, that is, of the parameters defining the position and orientation of the source and detector related to the rotational axis. However, when dealing with \replaced{experimental}{real} data obtained from \eg, a CT scanner, these parameters are provided a-priori and their accuracy is usually not enough to directly perform the reconstruction process, otherwise the obtained image will suffer from misalignment artifacts, impeding further analysis of the results. Therefore, an accurate (and computationally feasible) estimation of these parameters, or at least of a subset of them large enough to provide some desired reconstruction quality, is indispensable. 

The geometry of one cone-beam acquisition \azul{within the configuration described above} can be fully described by a set of seven parameters as follows. The 3D position of the center of the detector (where the central ray or focal axis is supposed to intersect), the 3D normal direction of the 2D plane containing the detector and the distance source-rotational axis, where the object will be placed. In this work, we propose an automatic, low-cost (computationally) and easy-to-implement methodology based on the variable projection optimization approach to estimate the detector shift in the direction orthogonal to both the focal and rotational axes, that is in the $u-$direction in Figure \ref{fig.conebeam}, and the in-plane detector rotation (over the focal axis) $\eta$ in Figure \ref{fig.conebeam}. These two parameters are well-known to be the cause of artifacts in reconstructions and thus affecting image quality and resolution, while the others are important if the goal is to perform metrological studies, as pointed out in \cite{cirp, lesaintcalibration} and the references therein. These claims are supported by numerical simulations in the cited works, however we can add some further explanations. It is easy to observe that a detector shift on $z$ and an inaccurate distance source-rotational axis will only have an effect on the scale of the reconstruction. The out-of-plane rotations of the detector will result in topological transformations of the reconstructed volume, with negligible effects on \eg\ number of features, as well as a detector shift on $y$. Note that we will not support these observations either analytically or numerically. \replaced{On the other hand, either}{Finally, both} a detector shift on $x$ \replaced{or}{and} a detector in-plane rotation will actually \emph{misalign} the rotational axis \azul{with respect} to the center of sinograms, producing double-edge artifacts as it will be observed in our experiments. \azul{Finally, the tomographic rotation angle is the remaining parameter needed to fully characterize a projection in a CT experiment. An inaccurate knowledge of it can also be a source of artifacts, this problem is addressed in \eg \cite{bresler}.}

\paragraph*{Related work \deleted{, contributions and overview}}
The estimation of the mentioned seven cone-beam geometrical parameters, although is possible, requires \eg\ non-automatic methods based on designing a calibration object usually composed of dense spheres with known geometry and relative location, and then, matching analytical projections with the acquired images \cite{massi, rishengict}. The obvious drawbacks of these methods are the need to design and manufacture the object and more importantly, the need to measure it often enough to account for non-static errors while non-reproducible errors are impossible to estimate.

Within the category of automatic methods (no \azul{pre-scan of a} reference object needed), a common approach is to iteratively perform the tomographic reconstruction of the object (or some relevant 2D slices) and update the geometric parameters until some criteria is met. A straightforward approach is to define a quality metric of the reconstruction and optimize it with some iterative approach, as done in \cite{variance} where the image contrast is considered, in \cite{entropy} its entropy or in \cite{l2gradient} the $l_2$ norm of its gradient. Following the same principle, in \cite{projectionmatch} a distance between the forward projection of the reconstruction and the measured data is computed and minimized with a quasi-Newtom algorithm. This reconstruction-based class of algorithms, although they are robust and obtain accurate results, have the drawback of being computationally expensive, as at every iteration a reconstruction (and eventually a forward projection) needs to be computed which is well-known to be costly mainly because of the \emph{backprojection} operation and specially when dealing with high-resolution data. The backprojection operation is the basis for almost any reconstruction method like the \emph{filtered backprojection} (FBP) algorithm \cite{kak} with a computational cost of $O(N^3)$ for $N^2$ detecting pixels and $N$ tomographic rotations to backproject a single 2D slice \cite{siamfast}. In addition, the  global convergence rate may also be slow as these algorithms only depend on the defined loss-function and its gradient (if available) and not on any geometrical prior information of the imaging system.  

Another class of methods not relying on reference objects are based on the so-called consistency conditions of the imaging geometry. They exploit relationships between noiseless and perfectly aligned projections that characterize the range of the imaging operator. In parallel-beam tomography, the alignment problem is usually reduced to the \emph{center of rotation} estimation. Necessary and sufficient consistency conditions to characterize the range of the Radon transform are well studied and known as the Helgason-Ludwig conditions \cite{helgason}. They are very useful for the center of rotation estimation problem and state for parallel projections $p(\cdot,\theta)$ at view angle $\theta$ that the $n-$moment $\int_\R p(t,\theta)t^n \, dt$ is a homogeneous polynomial of degree $n$ in $\{\cos\theta,\sin\theta\}$. Yet another necessary but not sufficient condition is the symmetry relationship 
\begin{equation}
    p(t,\theta) = p(-t,\theta+\pi), \quad \text{for all} \ t\in\R, \theta\in[0,2\pi),
    \label{eq.parallelsym}
\end{equation}
where $\theta$ is to be understood modulo $2\pi$ as $p$ is $2\pi-$periodic in $\theta$. In such geometry, based on (\ref{eq.parallelsym}), the estimation of the center of rotation shift  is easy to implement based on cross-correlation by doing 
\begin{equation*}
\textstyle{\frac{1}{2}}\argmax\limits_{t\in\R} \left\{ p (\cdot, 0) \star p (-\cdot, \pi) \right\}.
\end{equation*}  

In addition, the computing is low-expensive as it uses the fast Fourier transform to perform the cross-correlation operation $\star$.
In fan-beam geometries, the analogue of the Helgason-Ludwi\replaced{g}{n} consistency conditions were presented in \cite{finch} but no dedicated optimization strategy has been proposed to our knowledge. A symmetry relationship for fan-beam projections is also available and will be written out below in (\ref{eq.fancondition}) as it will be the basis for our fan-beam sub-problem. In \cite{yang2013}, the fan-beam center of rotation is estimated based on such relationship (\ref{eq.fancondition}), we will present an analysis of such method with a possible way to improve it. 

Cone-beam consistency conditions have been proposed in the last decade, in \cite{clackdoyle} necessary conditions are presented as homogeneous polynomials in $\{\cos\beta, \sin\beta\}$ with $\beta$ the rotation angle, in the same spirit of the Helgason-Ludwig conditions but without being sufficient as the Helgason-Ludwig's are. This result is with no doubts very relevant from a theoretical point of view but further analysis needs to be done with respect to the optimization strategy and algorithms to apply it in real tomographic experiments. In \cite{lesaintcalibration}, the 0th-order moment condition of \cite{finch} was used and extended to cone-beam projections by considering virtual detectors parallel to a line joining two cone-beam sources, from where a family of fan-beam projection are taken. Closely related, consistency conditions based on epipolar geometry were introduced in \cite{aichert} and further developed in \cite{risheng} obtained from Grangeat's theorem \cite{natteredwubbeling}. They need the computing of the derivative of epipolar plane integrals redundant for two sources. Both works \cite{lesaintcalibration, aichert} provide interesting tools to develop optimization algorithms, however they require extensive 3D interpolations and the use of projective geometry to compare every pair of 2D projections on the data. Such papers are mainly devoted to the description of the consistency conditions and how to compute them, however, adapted algorithms to estimate geometrical parameters are still to be investigated.

\paragraph*{\azul{Contributions and overview}} \azul{First, we propose two (to our knowledge) new, fast and automatic approaches to estimate the center of rotation in 2D fan-beam experiments based on sinogram registration and a fixed point iterative method. Both methods are based on the fan-beam symmetry relationship (\ref{eq.fancondition}) without the need of integrating over the angular variable such identity as is done in \cite{yang2013} and explained below. The remaining rebinning and signal registration strategies are addressed and detailed for both methods. Later, in cone-beam tomography,} we propose the use of the variable projection approach \cite{golub73} from continuous optimization. \azul{Variable projection solves structured optimization problems by \emph{removing} a subset of variables from the loss function which are relatively easy to solve and then it solves a \emph{reduced} problem over the remaining variables. This approach will be used} to estimate both shift $h$ and in-plane rotation $\eta$ parameters described above directly from \azul{cone-beam} projections. The inner problem (for $h$) will \azul{still} be based on \deleted{the fan-beam symmetry}relationship (\ref{eq.fancondition}) adapted to 3D \replaced{geometry}{data} \deleted{by studying and expanding Yang's method and by other two (to our knowledge) new approaches, based on 2D sinogram registration and a fixed point iterative method. The presented fan-beam work can already be used to estimate the center of rotation in 2D fan-beam experiments.} \azul{by expanding the corresponding rebinning and signal registration strategies of the two previously introduced 2D approaches.} Finally, \deleted{for cone-beam tomography,}the reduced problem on $\eta$ by \emph{projecting} $h$ will be solved via a gradient-descent algorithm with adaptive step size. The resulting approach will be proved to be less-expensive computationally with respect to state-of-the-art methods. Particularly, Section 2 is devoted to the fan-beam center of rotation alignment problem, \added{where the two strategies are formulated. They are} expanded to cone-beam \replaced{tomography}{data} in Section 3 \replaced{via}{followed by} the variable projection approach\deleted{to estimate $\eta$}. Section 4 presents numerical results \added{for both fan- and cone-beam geometries} with simulated and \replaced{experimental data from an industrial CT system}{industrial data} while conclusions are presented in Section 5. 

\section{2D problem: alignment in fan-beam tomography}

\newcommand{\quotient}{\R / 2\pi\mathbb{Z}}
\newcommand{\PIquotient}{\R / \pi\mathbb{Z}}

\noindent
This section studies the \emph{equidistant-detector fan-beam transform} $\fan\colon U \to V$ with circular source trajectory. It is defined through the \emph{divergent-beam transform} $D$ \cite{faridani} that acts on $f\in U$ as 
\begin{equation}
Df(a,\theta) = \int\limits_0^\infty f(a+t\theta) \, dt, 
\label{eq.divergentbeam}
\end{equation}
and gives the integral, or projection, of $f$ over the ray with direction $\theta\in\s^{1}$ and origin $a\in\R^2$. $\s^1$ is the unit circle, \ie, $\s^1 = \{ \theta\in\R^2, \| \theta\|_2=1 \}$. 

\begin{figure}[!t]
\centering
{\LARGE \scalebox{0.5}{     
\begin{tikzpicture}[x=1cm, y=1cm, scale=4]

        \draw [color=blue!20, fill=blue!20,rotate around={190:(0,0)}, scale=0.35] plot[smooth cycle, tension=1] coordinates {(1, 0) (0.6,0.4) (0.6,0.9) (0,1) 
        (-0.6,0.9) (-0.8,0.7) (-0.6,0.4) (-1,0) 
        (-0.8,-0.8) (0,-1) (0.8,-0.8)};

		
		\draw [>=latex,->, color=black] (-1.,0.) -- (1.5,0.);   
        \draw (1.5,0.) node[anchor=south east] {$x$};
		\draw [>=latex,->, color=black] (0.,-1.) -- (0.,2.);  
        \draw (0.,2.) node[anchor=north east] {$y$};

		
		\draw (0.,0.05) node[anchor=south west] {$f$};
				

        \def\x1{-1}
        \def\y1{-1.2}
        \def\ay2{1.5}

        \begin{scope}[rotate around = {25:(0,0)}]
        \foreach \t in {\x1,\x1+0.5,\x1+1,\x1+1.5,\x1+2} 
        { 
	    \draw[>=triangle 45,->, dashed, gray] (\x1+1,\y1) -- (\t,\ay2);
        }
        \draw (\x1+1,\y1) node[anchor=south west] {$r\beta$ (source)};
        \draw [very thick] (\x1-0.1,\ay2)-- (\x1+2.1,\ay2);	
        \draw [>=latex,->, thick] (0,0)--(0,-0.48) node[anchor=west] {$\beta$};

        \draw [>=latex,<-,color=black] (\x1-0.1,0)-- (\x1+2.1,0);  
        \draw [] (\x1-0.1,0) node[anchor=south] {$s$}; 
        \draw (\x1+0.9,\ay2) node[anchor=south] {\rotatebox{25}{$\fan f (\cdot,\beta)$}};
        \end{scope}		

\end{tikzpicture}}}
\caption{Fan-beam aligned tomographic acquisitions.}
\label{fig.fan}
\end{figure}
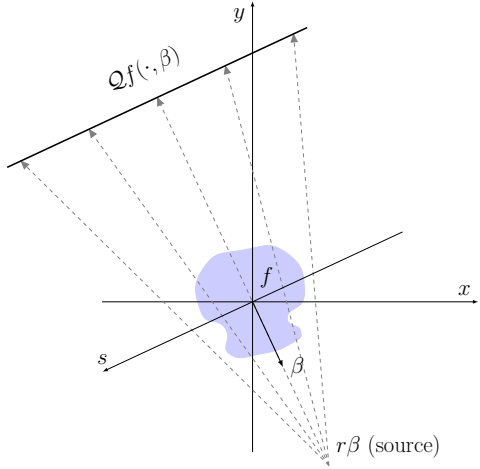

Particularly, if equidistant-detector fan-beam acquisitions are considered in (\ref{eq.divergentbeam}) as those illustrated in Figure \ref{fig.fan}, we parameterize with $(s,\beta) \in \R\times \s^1$ respectively the ray position over a linear detector (for clarity in notations, the $s - $axis crosses the origin and is parallel to the linear detector) and the tomographic rotation. Then, for $f \in U$, $\fan$ is defined as
\begin{equation*}
    \fan f(s, \beta) = Df(r \beta, \dfrac{s\beta^\perp-r\beta}{\|s\beta^\perp-r\beta\|_2}),   
\label{eq.fandefinition}
\end{equation*}
with 
$r>0$ denoting the fixed source-object distance, or source radius. That is,  $\fan f(s, \beta)$ gives the projection of $f$ following the ray starting at the source $r \beta$ towards the detecting point $s\beta^\perp$. $\beta^\perp$ is defined as the vector in $\s^1$ orthogonal to $\beta$ such that if $\beta=(\cos\varphi,\sin\varphi)$ then $\beta^\perp = (-\sin\varphi,\cos\varphi)$ with $\varphi$ the polar angle of $\beta$. The fan-beam transform can also be defined through the Radon transform and a rebinning operation as done in \cite{guerrero2023}.

To simplify notations, throughout the rest of the paper, $\beta$ will be parameterized by its polar angle which we will continue to denote by $\beta$.

The involved functional spaces are $U$, the \textit{feature} space that is the Schwartz space $\mathcal S(\R^2)$ of rapidly decreasing functions defined on $\R^2$; and $V$, the \textit{sinogram} space, is the Schwartz space for functions defined on $\R\times \s^1$.  



The range of $\fan$ will be denoted here by $\ran(\fan)\coloneqq\{\fan f\colon f\in U\}$. In the following, let $g\in\ran(\fan)$. We recall a useful symmetry property \cite{guerrero2023} of a sinogram $g$, that is, for any $(s,\beta)$ we have 
\begin{equation}
 g(s, \beta)  = g(-s,\beta+ \pi +2\arctan\dfrac{s}{r}). 
\label{eq.fancondition}
\end{equation}

Note that last condition is necessary but not sufficient to characterize $\ran(\fan)$. The fan-beam tomographic reconstruction problem is then posed as
\begin{equation*}
    \text{given} \ g\in \ran(\fan), \  \text{find} \ f\in U \  \text{such that} \ \fan f = g.
\end{equation*}

This problem will not be addressed in this work, and the reader can refer to \eg, \cite{guerrero2023, kak} for different techniques to solve it. Instead, the \emph{alignment} problem is introduced in the following.

Throughout the rest of this section, we will denote by $\g$ a \emph{misaligned} fan-beam sinogram in $s$, 
\ie, a sinogram translated in the form
\begin{equation}
     \tilde g(s,\beta) = g(s-h,\beta) \eqqcolon \tau_h g(s,\beta),  \quad \text{for all} \ (s,\beta), 
     \label{eq.shifted}
\end{equation}
for some $g\in \ran(\fan)$ and some \emph{unknown} $h\in\R$. $\tau_h$ is the translation operator by $h$ on the first variable. The alignment problem we address here is the estimation of such $h$.

In real scanning experiments, this happens when the center of rotation, expected to be at $s=0$ for any $\beta$, is shifted in the $s-$direction, and needs to be corrected before any reconstruction strategy. We can now state the problem as

\begin{pro}\label{pro.fan}
Given some misaligned data $\g$ in the form of (\ref{eq.shifted}), find $h\in\R$ such that
\begin{equation*}
\tau_h \g \in \ran(\fan).
\end{equation*}
\end{pro}




We propose to solve Problem \ref{pro.fan} based on relationship (\ref{eq.fancondition}). As this condition is not sufficient to characterize $\ran(\fan)$ we need to make some assumptions on the existence of the solution. Of course, if the data is simply in the form of (\ref{eq.shifted}), the condition is also sufficient, but this will never be the case with \replaced{experimental}{real} data.  

For any $(s,\beta)$, it is clear that, if the misalignment value $\h$ is known, from (\ref{eq.fancondition}) we have 
\begin{equation}
    \g(s+\h,\beta) = \g(-s+\h,\beta+\pi+2\arctan\dfrac{s}{r}), 
    \label{eq.shiftedins}
\end{equation}
then by calling $q = s+\h$, we also have
\begin{equation}
    \g(q,\beta) = \g(-q+2\h,\beta+\pi+2\arctan\dfrac{q-\h}{r}). 
    \label{eq.shiftedinq}
\end{equation}

The last equality will not only serve us in the derivation of the methods, but also to obtain a quality metric for any estimated $\h$, as we only need to interpolate on the right-side term and compare with the data $\g$ through some distance, as it will done below in the numerical results Section \ref{sec.2d.results}.


\subsection{Yang's method as a nearest neighbour interpolation} 

In \cite{yang2013}, Yang \ea. presented a method to solve Problem \ref{pro.fan} also based on (\ref{eq.fancondition}). Briefly, with different notations and explanations, the authors integrated over $\beta$ such identity and called 
\begin{equation}
    p(s) =  \int\limits_0^{2\pi} g(s, \beta) \, d\beta = \int\limits_0^{2\pi} g(-s,\beta+ \pi+2\arctan\dfrac{s}{r}) \, d\beta,
    \label{eq.yang_p}
\end{equation}
and claimed that 
\begin{equation}
   \int\limits_0^{2\pi} g(-s, \beta + \pi+2\arctan\dfrac{s}{r}) \, d\beta = 
   \int\limits_0^{2\pi} g(-s, \beta) \, d\beta = p(-s),
   \label{eq.yang_even}
\end{equation}
obtaining that $p$ should be an even function.

Yang's work suggested then to obtain $p$ by the first equality in (\ref{eq.yang_p}) as the sum of the fan-beam sinogram over $\beta$ for all $s$, and finally to estimate $\h$ by image registration between $p$ and the reversed signal $s\mapsto p(-s)$. 

To complement Yang's exposition, it is worth to note that identity (\ref{eq.yang_even}) is true because of the $2\pi-$periodicity of $g$ in $\beta$ and more importantly, only when $\beta \mapsto g(s,\beta)$ is defined on the whole $\s^1$ for any $s$. However, in real experiments we only dispose of a finite number of samples of $g$ and then (\ref{eq.yang_even}) is only an approximation, as follows. 

In the following, let $\gb\colon\Omega\to\R$ be a discrete fan-beam sinogram with domain the finite set of samples $\Omega=\{(s_i, \beta_j), i\in I, j\in J\}$ where $I, J$ index the detecting points and view angles respectively. Define the 1D array $\p  = \{\p_i\}$ by 
\begin{equation}
   \{\p_i \coloneqq \sum\limits_{j\in J} \gb(s_i, \beta_{j}), i \in I\}.
   \label{eq.p}
\end{equation}

In this discrete setting, 
Yang's method is actually based on approximation (\ref{eq.yang_even}) by
\begin{equation*}
   \p_i \approx \sum\limits_{j\in J} \gb(-s_i, \beta_{k_j}) = 
   \sum\limits_{j\in J} \gb(-s_i, \beta_j), \quad \text{for all} \ i \in I,
\end{equation*}
where $\beta_{k_j}$ is the nearest neighbour approximation of $\beta_j+ \pi +2\arctan\dfrac{s_i}{r}$ such that $(-s_i, \beta_{k_j}) \in \Omega$.

The last equality holds as for every $j\in J$, it corresponds a unique $k_j\in J$ (because of the angular $2\pi-$periodicity of $\gb$) and as the sum is carried out over all $j\in J$.
 
It is true that implementing such method is easy and the numerical computation is fast compared to other methods. Namely, we need to compute $\p$ as in (\ref{eq.p}), reverse its order and compute the shift between both 1D signals \eg, by cross-correlation-based signal registration \cite{registration}. However, the error related to the implicit nearest neighbour interpolation should be studied.

\subsection{Higher degree interpolation refinement} Nearest neighbour interpolation errors may be severe enough, specially when we only dispose of a sparse angular sampling of the data. Yang's method can easily be improved by higher degree interpolation as follows. Define the 1D array $\q  = \{ \q_i \}$ by
\begin{equation}
    \{\q_i \coloneqq \sum\limits_{j\in J} \gb(-s_i, \beta_{j} + \pi+2\arctan\dfrac{s_i}{r} ), i\in I \},
    \label{eq.q}
\end{equation}
where, for all $i\in I$, $\gb(s_i, \beta_{j} + \pi+2\arctan\dfrac{s_i}{r})$ is approximated in the $\beta-$variable by 1D interpolation of degree higher than $0$. We will show in the following that $\h$ can be estimated as the shift between $\p$ and $\q$ again by cross-correlation.

\added{Denote $\tilde p$ the analogue of Yang's 1D signal $p$ in (\ref{eq.yang_p}) obtained from a misaligned sinogram $\g$, \ie,}
\begin{equation*}
    \added{\tilde p(s) } \added{= \intpi \g(s, \beta) \, d\beta,}
\end{equation*}
\added{which, if $\h$ is known, verifies identity (\ref{eq.shiftedinq}) after $\beta-$integration (while $p$ verifies (\ref{eq.fancondition})). Thus, denote the 1D signal}
\begin{equation*}
    \added{w(s) = \intpi \g(-s,\beta+\pi+2\arctan\dfrac{s}{r}) \, d\beta.}
\end{equation*}

\added{Note that $\mathbf{w}$ is the discrete version of $w$, which is available from the data after interpolation. The following result holds.}
\deleted{Starting from (\ref{eq.shiftedins}), denote the operators}
\deleted{and the following result holds.}
\begin{prop}\label{pro.ly}
For $\g$ given as in (\ref{eq.shifted}), the \replaced{signals $\tilde p$ and $w$}{operators $\w_0 \g$ and $\z_0\g$} are related by 
\begin{equation*}
    \replaced{\tilde p}{\w_0\g}(s) = \replaced{w}{\z_0\g}(s-2\h), \quad \text{for all} \ s\in\R.
\end{equation*}
\end{prop}
\begin{proof}
\added{Similar to the analysis of Yang's method, the identity}
\azul{
\begin{equation}
     \intpi \g(-s+2\h, \beta+\pi +2\arctan\dfrac{s-\h}{r}) \ d\beta   =  \intpi \g(-s+2\h, \beta+\pi + 2\arctan\dfrac{s-2\h}{r}) \  d\beta
  \label{eq.linearyangproof}
\end{equation}
}\azul{is true because of the $2\pi-$periodicity of $\g$ in $\beta$ and because $\beta \mapsto \g(s,\beta)$ is defined on $\s^1$ for any $s$. By (\ref{eq.shiftedinq}), the left-side of (\ref{eq.linearyangproof}) is $\tilde p(s)$  and the right-side is $w(s-2\h)$ by definition, thus the results follows.}
\deleted{Let $q=s+\h$, from identity (\ref{eq.shiftedinq}) we have}
\deleted{therefore
$
\w_0\g(q) = \z_0\g(q-2\h)
$
for all $q\in\R$.}
\end{proof}
\deleted{As in the analysis of Yang's method, the last equality in (\ref{eq.linearyangproof}) holds because of the $2\pi-$periodicity of $\g$ in $\beta$ and because $\beta \mapsto \g(s,\beta)$ is defined on $\s^1$ for any $s$.}

\added{Proposition \ref{pro.ly} implies that the discrete signals $\p$ and $\mathbf{w}$ are \emph{shifted} by $2\h$.} We are now ready to state a first algorithm to estimate $\h$, referred to here as Linear Yang \added{(LY)}.
\vspace{15pt}
{\hrule width 0.4\textwidth}
\begin{alg}\replaced{LY for Problem \ref{pro.fan}}{Linear Yang} \\
Given fan-beam data $\gb$  
\begin{enumerate}[label=\footnotesize\arabic*), itemsep=-0pt]
    \item Compute $\p$ by (\ref{eq.p}) and $\q$ by (\ref{eq.q}), the latter by linear interpolation on $\beta$
    \item $\h = \textstyle{\frac{1}{2}}\,\text{shift}(\p,\q)$
\end{enumerate}
{\hrule width 0.4\textwidth}\vspace{10pt}
\end{alg}

Finally, the shift operation can be performed by cross-correlation \cite{registration} as
\begin{equation}
  \text{shift}(\p, \q) \coloneqq \argmax_{i\in I} \p \star \q,
  \label{eq.cross}
\end{equation}
with $\star$ the discrete one-dimensional cross-correlation operator. Note that sub-pixel results can by achieved by zero-padding in the frequency domain if the cross-correlation is computed with the Fourier Transform, as we suggest.

\subsection{Main approaches}

We observed that both Yang's original approach and its extension via linear interpolation have the drawback of not performing so well with low contrast data as we will show later in the numerical results. This is explained by the fact that these methods are based on an integral operation over the angular variable, therefore the details of each projection are \emph{averaged out} after the sum (\ref{eq.p}) (and (\ref{eq.q}) for \replaced{LY}{the Linear Yang method}). We propose in the following two methods that avoid this integral operation.

\subsubsection{2D sinogram registration} 

We present in this section a 2D sinogram registration method still based on identity (\ref{eq.shiftedins}) without the need of the integral over $\beta$. \added{Given a sinogram $\g$,} define now the \replaced{2D signal}{operator} 
\begin{equation*} 
\replaced{z}{\x_hg}(s,\beta) = \replaced{\g}{g}(-s\deleted{+h},\beta+\pi+2\arctan\dfrac{s}{r}),
\end{equation*}
then we have the following proposition. 
\begin{prop}\label{prop.2d}
    A sinogram $\g$ in the form of (\ref{eq.shifted}) and \replaced{$z$}{$\x_0\g$} are translated for all $(s,\beta)$ as
    \begin{equation}
        \g(s,\beta) \approx \replaced{z}{\x_0\g}(s-2\h,\beta- \dfrac{2\h}{r} ),
        \label{eq.approx2d}
    \end{equation}
with an error related to the first order approximation of the arc-tangent function in the $\beta-$variable of \replaced{z}{$\x_0\g$}.
\end{prop} 
\begin{proof}
    From identity (\ref{eq.shiftedinq}), we have 
    \begin{equation*}
        \g(\replaced{s}{q},\beta) = \g(-\replaced{s}{q}+2\h,\beta + \pi + 2\arctan\dfrac{\replaced{s}{q}-\h}{r} ) \approx \g(-\replaced{s}{q}+2\h,\beta + \pi + 2\arctan\dfrac{\replaced{s}{q}}{r}-\dfrac{2\h}{r}),
    \end{equation*}
    where approximation (\ref{eq.approx2d}) yields for all $(\replaced{s}{q},\beta) \in \R\times\s^1$ after the use of the first order approximation of the arc-tangent function in the $\beta-$variable. Note that we always have $\left| s \right| / \replaced{r}{d} < 1$ and that there is no approximation in the $s-$variable.
\end{proof}

For the purposes of the present approach, even if the approximate translation in $\beta$ also gives us an estimation of $\h$, we will recover $\h$ from the $2\h$ translation value in $s$ of (\ref{eq.approx2d}) as this translation is exact. In the discrete setting detailed above, define \azul{the discrete version of $z$, \ie\ }the 2D array $\sd  = \{ \sd_{i,j} \}$ by
\begin{equation}
    \{ \sd_{i,j} \coloneqq \gb(-s_i, \beta_j + \pi+2\arctan\dfrac{s_i}{r} ), i\in I, j\in J \},
    \label{eq.s}
\end{equation}
computed by some interpolation method on the angular variable. The following  algorithm summarises the method, referred to as 2D registration \added{(2DR)}.
\vspace{12pt}
{\hrule width 0.4\textwidth}
\begin{alg}\replaced{2DR for Problem \ref{pro.fan}}{2D registration}\label{alg.2d} \\
Given fan-beam data $\gb$  
\begin{enumerate}[label=\footnotesize\arabic*), itemsep=-0pt]
    \item Compute $\sd$ as in (\ref{eq.s}) by (linear) interpolation on $\beta$
    \item $\h = \textstyle{\frac{1}{2}}\shift\limits_{s}(\gb,\sd)$
\end{enumerate}
{\hrule width 0.4\textwidth}\vspace{10pt}
\end{alg}

Here, the operation $\shift\limits_{s}$ returns only the shift value of the first variable $s$ of the 2D arrays, it can be computed by cross-correlation in the same way as (\ref{eq.cross}) with 2D signals ignoring the shift value in $\beta$. Zero-padding the frequency variable for sub-pixel accuracy is also possible with some computing memory requirements. Reduced memory algorithms for 2D image registration by cross-correlation are also available and presented in \cite{manuel}. Finally note that the contrast of each projection is taken into account in the $\shift\limits_{s}$ operation, resulting in a more robust method for \eg\ low-contrast data.


\subsubsection{A fixed point method}\label{sec.fixedpoint}

We propose in this section to fix $\beta$ \eg, $\beta=0$ in (\ref{eq.shiftedins}) and (\ref{eq.shiftedinq}) and analyse the translation properties between the resulting 1D signals. Let us denote the operators 
\begin{equation} 
\Lambda g(q) = g(q, 0), \quad \Pi_hg(q) = g(-q+2h,\pi+2\arctan\dfrac{q-h}{r}).
\label{eq.operators}
\end{equation}

\azul{Note that by (\ref{eq.shiftedinq}), the searched shift value $\h$ verifies 
$\Lambda \g =  \Pi_{\h} \g$.} We could already use a brute force approach and consider the least-squares problem 
\begin{equation*}
\min\limits_h \big\{\,L(h) \coloneqq \| \Lambda \g - \Pi_h \g \|_2^2 \, \big\},
\label{eq.leastsquares}
\end{equation*}
and use a gradient-based algorithm to solve it. However, with such an approach we are not taking into account any geometrical properties of $\Lambda \g$ and $\Pi_h \g$, and the method will rely solely on properties of the gradient of $L$. 

Instead, we will use the fact that $\Lambda\g$ is an approximate translated version of $\Pi_0\g$ by $2\h$ and analyse the error. We can then estimate the shift value again by signal registration procedures followed by an iterative refinement, as follows. 

\begin{prop} For a misaligned sinogram $\g$ and for all $s\in\R$, we have the approximate translation
\begin{equation}
   \Lambda\g(s) \approx \Pi_0\g(s-2\h),
   \label{eq.shiftfp}
   \end{equation}
with an $L_\infty$ error bounded by
\begin{equation*}
    \max\limits_{s} \left| \Lambda\g(s) - \Pi_0\g(s-2\h) \right| \leq C_{\h},
\end{equation*}
with
\begin{equation}
    C_{\h}  = \max\limits_{s,\beta}  \left| \g(s,\beta) - \g(s,\beta+ \dfrac{2\h}{r}) \right|.
    \label{eq.bound}
\end{equation}
\end{prop}

\begin{proof}
    From (\ref{eq.shiftedinq}), we have
\begin{equation}
    \g(q,0) = \g(-q+2\h, \pi +  2\arctan\dfrac{q-\h}{r}) \approx \g(-q+2\h, \pi+2\arctan\dfrac{q-2\h}{r}), 
  \label{eq.approx}  
\end{equation}
therefore $\Lambda\g(q) \approx \Pi_0\g(q-2\h)$ for all $q\in\R$. \azul{Note that the last approximation simply reads $\Pi_{\h}\g(q) \approx \Pi_0\g(q-2\h)$, which is the basis of the presented method as the dependency of $\h$ is not anymore in the operator but in the argument of the signal resulting in a \emph{shifted} signal with respect to $\Lambda \g$.}

For the error bound, the approximation in (\ref{eq.approx}) is in the argument of the arc-tangent function \azul{by $\h/r$}. It is very similar to the manipulation in (\ref{eq.linearyangproof}) that turns out to be an equality after the integral operation. Recall that $q=s+\h$, the error in $\beta$ is $2 \left| \arctan s/r-\arctan (s-\h) /r \right|$, and it can be shown to be bounded by $2\left| \h \right| / r$ after simple manipulations using the subtraction property of the arc-tangent function \cite[Equation~4.4.34]{abramowitz}, from where we obtain the bound (\ref{eq.bound}). Note finally that $\left| \h \right|/r \ll 1$ and that the maximum is attained in (\ref{eq.bound}) because $\g$ belongs to a Schwartz space. 
\end{proof}

We are able to provide a first rough approximation of $\h$ as 
\begin{equation}
      \h \approx\textstyle{\frac{1}{2}} \, \text{shift}(\Lambda\g, \Pi_0\g) \coloneqq 
      \textstyle{\frac{1}{2}} \argmax\limits_{s} \Lambda\g \star \Pi_0\g,
      \label{eq.crosscont}
\end{equation}
where here the shift and $\star$ operations are to be understood for signals defined on $\R$ and not discrete signals as in (\ref{eq.cross}), we keep the same notation for sake of simplicity.

Approximation (\ref{eq.crosscont}) is much worse than the estimation based on (\ref{eq.approx2d}) in \replaced{2DR}{the 2D registration method} as such translation property is exact in $s$ using the whole 2D sinogram. The translation error in $\beta$ in such approximation is reflected directly in the $s-$variable in (\ref{eq.shiftfp}) as $\{\Lambda\g,\Pi_0\g\}$ are 1D signals. We propose to refine (\ref{eq.crosscont}) iteratively. However, we need to make some assumptions on the error $C_{\h}$ in (\ref{eq.bound}). 

Denote the real function
\begin{equation*}
  T_g(h) = h + \shift(\Lambda g, \Pi_h g).
\end{equation*} 

From (\ref{eq.shiftedinq}), we observe \azul{again} that $\shift(\Lambda\g,\Pi_{\h}\g) = 0$, thus $\h$ is a solution of 
\begin{equation*}
    T_{\g}(h)=h,
    \label{eq.fixedpoint}
\end{equation*}
\ie, $\h$ is a fixed point of $T_{\g}$.

 We make the assumption for the error $C_{\h}$ that is such that $T_{\g}$ is a contraction in a neighbourhood of $\h$. 
 With such assumption, we can now present our fixed point method with trivial proof using the contraction theorem \cite{kantorovich}. 
\begin{prop}\label{prop.fixedpoint}
Problem \ref{pro.fan} has a unique solution in a neighbourhood of $\h$ given by the limit of the sequence $(h_k)$ defined by the iteration
\begin{equation}
h_0 = 0, \quad h_{k+1} = T_{\g}(h_k).
\label{eq.iterativefixedpoint}
\end{equation}
\end{prop}

At every call of the function $T_{\g}$ we need to compute $\{\Lambda\g,\Pi_h\g\}$ defined in (\ref{eq.operators}) which are 1D signals that after discretization can be computed by (linear) interpolation. An additional advantage of expressing them in terms of $q$ and not $s$ is that $\Lambda$ no longer depends on $h$ and no interpolation is needed, then it is pre-computed before the loop in the following algorithm\added{, referred to as Fixed point (FP).}
\newpage
{\hrule width 0.4\textwidth}
\begin{alg}\replaced{FP for Problem \ref{pro.fan}}{Fixed point}\label{alg.fp} \\
Given fan-beam data $\gb$ and $h_0=0$  
\begin{enumerate}[label=\footnotesize\arabic*), itemsep=-0pt]
    \item Set $\Lambda_i = \gb(s_i,0)$ for all $i\in I$ 
    \item Until convergence, for $k=0,1,\dots,$ do
    \begin{enumerate}[label=\footnotesize\arabic*), itemsep=-0pt]\vspace{-2pt}
    \item Compute $\Pi_i  = \gb(-s_i+2h_k,\pi+2\arctan\dfrac{s_i-h_k}{r})$ by (linear) interpolation for all $i\in I$
    \item $h_{k+1} = h_k + \textstyle{\frac{1}{2}}\shift(\Lambda, \Pi)$
    \end{enumerate}
    \item $\h = h_{k+1}$
\end{enumerate}
{\hrule width 0.4\textwidth}\vspace{5pt}
\end{alg}

\added{Finally, we can improve the robustness of FP by simply taking the median of $K$ calls of FP with different values of $\beta$ uniformly distributed in the angular range, then an eventual defective projection at $\beta=0$ will have no effect. Recall that $\beta=0$ is fixed in the derivation of FP. This modified algorithm is referred to as FP$_{K}$.}

\subsection{Discussion and computational cost}


We make the assumption that there exists an optimal value $\h$ as a solution to Problem \ref{pro.fan}. Yang and \replaced{LY}{Linear Yang (LY)} methods use the sinogram \emph{averaged} in $\beta$ in (\ref{eq.p}) to compute the shift, then both methods assume that the solution is findable after averaging the data. This is not needed for \replaced{2DR}{2D registration (2DR)} that performs the registration directly on the sinogram and \replaced{FP}{Fixed point (FP)} that does it on 1D projections. Note that FP also needs the additional contraction assumption of Proposition \ref{prop.fixedpoint} requiring the first iteration to register the approximately shifted projections with a value \emph{close enough} to $\h$. \added{These observations suggest that Yang and LY methods are robust to noise, thanks to the averaging operation, as well as 2DR, because of the 2D registration on the full sinogram. On the other hand, FP is more sensitive to noise, but can be easily improved with its extended version FP$_K$. Also, Yang and LY need the solution to be findable after the averaging operation, which makes them less attractive for challenging data, \eg\ low-contrast data or data with few details. FP and 2DR operate directly on the raw data, therefore \emph{any} detail will be taken into account by the algorithm. These observations will be confirmed in the numerical results Section \ref{sec.results} below.}

The complexity analysis is done for a sinogram with $N$ pixels and $N_\beta$ rotations. 1D signal registration is computed with complexity $O(N\log N)$ if the Fast Fourier Transform (FFT) is used, then this is the complexity of Yang's method. LY  has a complexity $O(NN_\beta)$ as this is the cost of the interpolation in (\ref{eq.q}). 2DR also has a complexity $O(NN_\beta)$ in the interpolation (\ref{eq.s}) but $O(NN_\beta\log(NN_\beta))$ for the 2D signal registration step again with FFT. Finally, FP with $M$ iterations has a complexity $O(MN\log N)$ with FFT at every iteration for the registration. Note that $M=5\ll N_\beta$ in all our numerical experiments. If $M \ll N_\beta$, FP is less expensive computationally than 2DR, and if $M$ is bounded, FP has the same low complexity of Yang's method. \added{FP$_{K}$ still has lower cost than 2DR as long as $KM\ll N_\beta$}. In any case, any strategy based on iterative reconstructions has a cost $O(MN^2N_\beta)$ for $M$ iterations, due to the backprojection step. Although there exist techniques to accelerate the backprojection \cite{guerrero2023, george2013fast}, all the methods presented here are drastically faster.


\section{3D \replaced{problem}{case}: alignment in cone-beam tomography}
\noindent
In this section the \emph{cone-beam transform} with circular source \azul{trajectory} will be studied, it is defined again by means of the divergent-beam transform (\ref{eq.divergentbeam}) with source $a\in\R^3$ and ray direction $\theta\in\s^2$. It is parameterized by 2D detector coordinates $(u,v)\in\R^2$ and tomographic angle $\beta\in\s^1$ (parameterized again by its polar angle as the fan-beam case) as illustrated in Figure \ref{fig.conebeam}. It is denoted and defined for $f\in U$ as 
\begin{equation*}
    \cone f(u,v,\beta) = Df(ra_\beta,  \theta_{u,v,\beta}),
\end{equation*}
where $a_\beta = (\cos\beta,0,\sin\beta)$ and $\theta_{u,v,\beta} = \dfrac{ \bm u  - ra_\beta  }{  \| \bm u - ra_\beta \|_2 } $, $\bm u$ is the 3D Cartesian coordinates of the detecting point $(u,v)$ and rotation $\beta$, see \cite[Section~5.5.1]{natteredwubbeling} for the expression of $\bm u$ in terms of $u,v$ and $\beta$.

\paragraph*{Notations} In the following, $g\in\ran(\cone)$; the translation (in $u$) and rotation operators will be denoted by
\begin{equation*}
    \tau_{h} g(u,v,\beta) = g(u-h, v,\beta), \quad
    \kappa_\eta g(u, v,\beta) = g(u\cos\eta-v\sin\eta, u\sin\eta+v\cos\eta,\beta). 
\end{equation*}

Note that operator $\kappa_\eta$ is given by multiplying the rotation matrix 
$
\begin{pmatrix}
    \cos\eta & -\sin\eta \\
    \sin\eta & \cos\eta
\end{pmatrix}
$
to the detector coordinates $(u,v)$ of $g$ for all $\beta$. The misaligned cone-beam projections considered here will be denoted by $\g$ which are obtained from $g\in\ran(\cone)$ in the form
\begin{equation}
    \g = \tau_h\kappa_\eta g,
    \label{eq.miscone}
\end{equation}
for some unknown $(h,\eta)$. That is, we consider cone-beam projections measured with a relative rotation between detector and rotation axis, followed by a shift of the detector in the $u-$direction. We can now state the cone-beam alignment problem as the following.

\begin{pro}\label{pro.cone}
Given misaligned cone-beam data $\g$ in the form of (\ref{eq.miscone}), find $h\in\R$ and $\eta \in \left(-\pi/2,\pi/2\right)$ such that
\begin{equation*}
\kappa_\eta \tau_h \g \in \ran(\cone).
\end{equation*}
\end{pro}

We will make use of the fan-beam symmetry relationship (\ref{eq.fancondition}) extended to cone-beam projections with $v=0$, \ie, for every $(u,\beta)$ we have
\begin{equation*}
    g(u, 0, \beta)  = g(-u,0,\beta + \pi + 2\arctan\dfrac{u}{r}), 
\end{equation*}
which can be expressed in terms of $\g$ and the correct alignment values $(\h,\eta^*)$ by noting that $ \kappa_{\eta}^{-1} \tau_{h}^{-1}  \g  =  \kappa_{-\eta} \tau_{-h}  \g = g $  and then
\begin{equation}
   \kappa_{-\et} \tau_{-\h}  \g(u,0,\beta)  =   \kappa_{-\et} \tau_{-\h}  \g(-u,0,\beta + \pi + 2\arctan\dfrac{u}{r}).
   \label{eq.cone.sym_tildeg}
\end{equation}


\subsection{A variable projection method}

We propose in this section to pose Problem \ref{pro.cone} as a joint least-squares problem in $(h,\eta)$ and solve it by the variable projection (VP) approach \cite{golub73} \azul{as follows.} Denote now the operators
\begin{equation}
\begin{array}{ll}
\Lambda_{h,\eta}g(u,\beta) \!\!\!\!\! &= \kappa_{-\eta} \tau_{-h}  g(u,0,\beta), \\
\Pi_{h,\eta}g(u,\beta) \!\!\!\!\! &= \kappa_{-\eta} \tau_{-h}  g(-u,0,\beta + \pi + 2\arctan\dfrac{u}{r}), 
\end{array}
\label{eq.cone.operators1}
\end{equation}
\azul{where, by (\ref{eq.cone.sym_tildeg}), we have $\Lambda_{\h,\et}\g = \Pi_{\h,\et}\g$.} Then\deleted{, by (\ref{eq.cone.sym_tildeg}),} $(h,\eta)$ can be estimated by solving 
\begin{equation}
\min_{h,\eta}\big\{ L(h,\eta) \coloneqq \|  \Lambda_{h,\eta}\g - \Pi_{h,\eta}\g \|^2_2  \big\},
\label{eq.cone.leastsquares}
\end{equation}
which again can be solved with any gradient-based algorithm. However, problem (\ref{eq.cone.leastsquares}) can be severely ill-posed due to the non-guarantee of convexity of $L$ and to differences in sensitivity of $h$ and $\eta$, as pointed out in \cite{tristan2018} in a different reconstruction problem. In addition, this brute-force strategy won't allow us to use results of the 2D fan-beam case as the VP approach will, as presented in the following. 

Inspired by \cite{tristan2018}, which uses the VP approach to solve a joint and generalized tomographic alignment-reconstruction problem, we will project $h$ onto $\eta$ by setting 
\begin{equation}
    \bar h(\eta) = \argmin_h L(h, \eta),
    \label{eq.cone.fanprojected}
\end{equation}
which is essentially a \emph{tilted} fan-beam alignment problem that can be solved based on the approaches presented in the previous section with the incorporation of the tilting variable $\eta$ on the fan-beam geometry as done below. Then, the VP approach states that problem (\ref{eq.cone.leastsquares}) is equivalent to the reduced problem
\begin{equation}
    \min_{\eta} \big\{ \bar L(\eta) \coloneqq  L(\bar h(\eta),\eta) \big\}.  
    \label{eq.cone.reduced}
\end{equation} 

The reduced problem (\ref{eq.cone.reduced}), expected to be less ill-conditioned, as for the reconstruction problem in \cite{tristan2018}, can already be solved by gradient-based methods. Indeed, we have the following expression of the derivative of $\bar L$, denoted $\nabla \bar L$, obtained from \cite{2008gradient, tristan2012},
\begin{equation*}
    \nabla \bar L(\eta) = \nabla_\eta L(\bar h(\eta), \eta).
\end{equation*}

We make the assumption that optimal values $(\h,\et)$ exist for Problem \ref{pro.cone} that minimizes $L$ and that $L$ is convex in a neighbourhood of $(\h,\et)$. Following \cite{tristan2018}, it is guaranteed that $\bar L$ is also locally convex and that a local minimum $\eta^*$ of $\bar L$ with the corresponding $\h = \bar h(\eta^*)$ is a local minimum of $L$. The second derivative of $\bar L$ is also available in \cite{2008gradient, tristan2012} that would allow us to use second-order methods to solve the reduced problem (\ref{eq.cone.reduced}), however, as the reduced problem is univariate, we limit ourselves to first-order methods.

\paragraph*{The VP iteration}

We will adapt \replaced{2DR}{the 2D registration} (Algorithm \ref{alg.2d}) and \replaced{FP}{fixed point} (Algorithm \ref{alg.fp}) fan-beam approaches to solve the sub-problem (\ref{eq.cone.fanprojected}), to do this we need to incorporate the tilting variable $\eta$ in the fan-beam geometry. Let us rewrite the operators defined in (\ref{eq.cone.operators1}) after setting $q=u+h$ as
\begin{equation}
\begin{array}{rl}
   \Lambda_{\eta}g(q,\beta)  = \!\!\!\! &   g(q\cos\eta, -q\sin\eta,\beta), \\
   \Pi_{h,\eta}g(q,\beta)  = \!\!\!\! &  g((-q+2h)\cos\eta, (q-2h)\sin\eta, \beta+ \pi + 2\arctan\dfrac{q-h}{r}).
\end{array}
\label{eq.cone.operators2}
\end{equation}

Note that $\Lambda_\eta$ is not parameterized with $h$ anymore and we can observe a shift in the $q-$variable by $2\h$ between both signals because of (\ref{eq.cone.sym_tildeg}). Then, we have the following extension of Proposition \ref{prop.2d}, with similar proof up to easy manipulations.
\begin{prop}\label{prop.cone.2d}
Let $\eta\in\s^1$ be fixed, the 2D signals $\Lambda_{\eta}\g$ and $\Pi_{0,\eta}\g$ are translated one to another for all $(u,\beta)\in\R\times\s^1$ as
    \begin{equation*}
        \Lambda_\eta\g(u,\beta) \approx \Pi_{0,\eta}\g(u-2\h,\beta- \dfrac{2\h}{r} ),
        \label{eq.2dshift}
    \end{equation*}
with an error in the $\beta-$variable of the same order of approximation (\ref{eq.approx2d}).
\end{prop}

\replaced{FP}{The fixed point method} is easily extendable to a tilted fan-beam as well. Indeed, from operators (\ref{eq.cone.operators2}) and setting $\beta=0$ we have the following result.
\begin{prop}\label{prop.cone.fp} 
For a fixed $\eta\in\s^1$, we have for all $u\in\R$, 
\begin{equation*}
   \Lambda_\eta\g(u,0) \approx \Pi_{0,\eta}\g(u-2\h,0),
   \end{equation*}
with an $L_\infty$ error bounded as in (\ref{eq.bound}).
\end{prop}

Therefore, we can estimate $\bar h(\eta)$ for a fixed $\eta$ in (\ref{eq.cone.fanprojected}) by refining iteratively 
\begin{equation*}
   \textstyle\frac{1}{2} \shift(\Lambda_\eta\g(\cdot,0), \Pi_{0,\eta}\g(\cdot,0))
\end{equation*} 
following the same procedure as done in the \replaced{derivation of FP}{the fixed point method} of Section \ref{sec.fixedpoint}. That is, we can apply Algorithms \ref{alg.2d} and \ref{alg.fp} to solve the inner problem (\ref{eq.cone.fanprojected}) of the VP strategy based respectively on Propositions \ref{prop.cone.2d} and \ref{prop.cone.fp}. Particularly and for a fixed $\eta$, the 2D fan-beam data $\{\gb,\mathbf{z}\}$ in Algorithm \ref{alg.2d} are the discrete arrays
\begin{equation*}
\gb =\{\Lambda_\eta\g(u_i,\beta_j)\}, \quad \mathbf{z} = \{\Pi_{0,\eta} \g(u_i,\beta_j)\}, \quad i\in I, j\in J,     
\end{equation*}
with $I,J$ indexing the detector columns and view angles respectively. Whereas in Algorithm \ref{alg.fp}, the 1D signals $\{\Lambda,\Pi\}$ are
\begin{equation*}
\Lambda =\{\Lambda_\eta\g(u_i,0)\}, \quad \Pi = \{\Pi_{0,\eta} \g(u_i,0)\}, \quad i\in I.    
\end{equation*}

We are ready to state the VP method suggested to solve Problem \ref{pro.cone} by solving both the inner problem (\ref{eq.cone.fanprojected}) and the reduced problem (\ref{eq.cone.reduced}). We will denote by $\textsf{h}(\gb,\eta)$ the \azul{computational} strategy to solve the inner problem for some cone-beam data $\gb$ and a fixed $\eta$ that provides an estimate of the shift value $\bar h(\eta)$. It can be done either via the 3D-adapted Algorithms \ref{alg.2d} or \ref{alg.fp} as described above. The \emph{univariate} reduced problem can be solved via \eg, gradient-descent as in the following algorithm\added{, referred to VP}. 

\vspace{15pt}
{\hrule width 0.4\textwidth}
\begin{alg}\replaced{VP for Problem \ref{pro.cone}}{Variable projection}\label{alg.vp} \\
Given cone-beam data $\gb$ and $\eta_0=0$   
\begin{enumerate}[label=\footnotesize\arabic*), itemsep=-0pt]
    \item Until convergence, for $k=0,1,\dots,$ do
    \begin{enumerate}[label=\footnotesize\arabic*), itemsep=-0pt]\vspace{-0pt}
    \item $h_k = \textsf{h}(\gb, \eta_k)$  
    \item $\eta_{k+1} = \eta_{k} - \gamma_k \nabla_\eta L(h_k, \eta_k) $, for some $\gamma_k > 0$ 
    \end{enumerate}
    \item $(\h,\et) = (\textsf{h}(\gb, \eta_{k+1}), \eta_{k+1}$)
\end{enumerate}
{\hrule width 0.4\textwidth}
\end{alg}

\subsection{Discussion}

\paragraph*{Gradient descent and step sizes $\gamma_k$} We can write down the derivative $\nabla \bar L(\eta) = \nabla_\eta L(\bar h(\eta), \eta)$ simply by applying the multivariable \emph{chain rule} as
$
\nabla \bar L (\eta) = 2\langle G(\eta), \nabla G (\eta)\rangle, 
$
with the function $G$ and its gradient expressed respectively as 
\begin{equation*}
\begin{split}
    G(\eta) \coloneqq & \Lambda_{\eta}g(q,\beta) - \Pi_{h,\eta}g(q,\beta), 
    \\
    \nabla G (\eta) = & \Lambda_{\eta} \left(-q(\sin\eta,\cos\eta,0) \cdot \nabla g  \right)(q,\beta) - \Pi_{h,\eta} \left((q-2\bar h(\eta))(\sin\eta,\cos\eta,0) \cdot \nabla g\right) (q,\beta), 
\end{split}
\end{equation*}
where $\nabla g = (\nabla_u g, \nabla_v g, \nabla_\beta g)$ is the gradient of the data, to be calculated \eg, by finite differences. However, with high-resolution data care is to be taken due to computing (and storing) $\nabla g$ can be not efficient computationally, as well as the interpolations needed to compute $\nabla G$. As an alternative, finite differences directly on $\bar L$ is possible with care on choosing the step size $\Delta\eta$ depending on the resolution of the data. The latter approach is the one we used in our experiments with $\Delta\eta=0.001$ radians after verification of not oversampling while still maintaining a desired accuracy.     

The step sizes $\gamma_k$ are selected with a line search approach following a sufficient decrease and backtracking strategy until the Armijo condition is verified \cite[Algorithm~3.1]{nocedal}. Namely, if $\gamma_k$ does not verify the Armijo condition, we decrease $\gamma_k$ by a contraction factor of $1/2$. Let us recall that this step size selection produce showed crucial importance in the convergence of the gradient descent algorithm in our experiments, in accordance with optimization theory. A fixed $\gamma$ or even $\gamma_k$ producing insufficient decrease of $\bar L$ at each iteration, could translate to a non-convergent algorithm.

\section{\deleted{Cone-beam} Numerical results}\label{sec.results}
\noindent
Validation of the methods is done with simulated and experimental data \added{for both fan- and cone-beam geometries. All the presented algorithms are tested and compared with an additional non-automatic method for the cone-beam case.} 

\subsection{\added{Fan-beam data}}\label{sec.2d.results}
\subsubsection{\added{Simulated experiment}}
\begin{figure}[!t]
    \centering
    \includegraphics[width=0.155\textwidth]{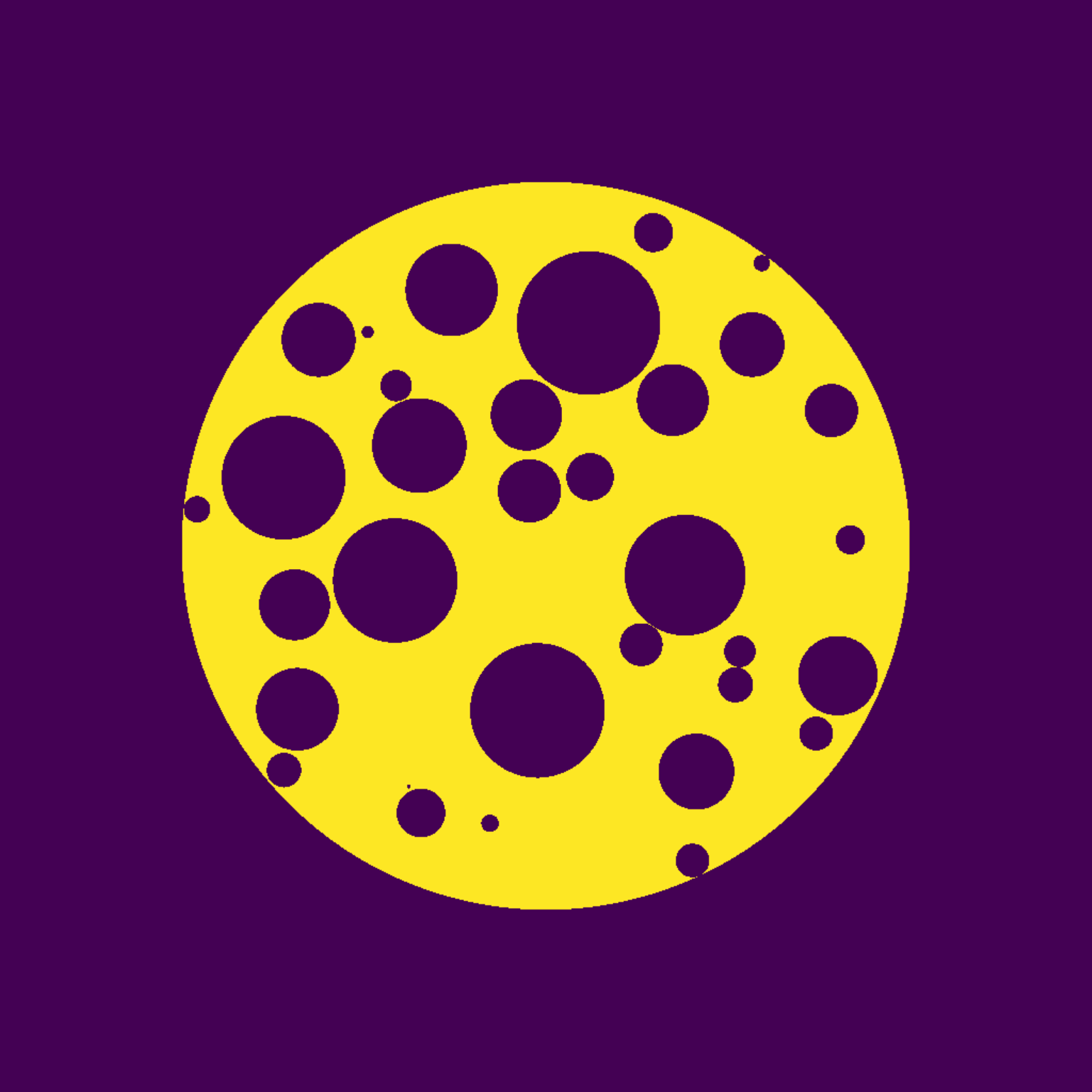}
    \includegraphics[width=0.155\textwidth]{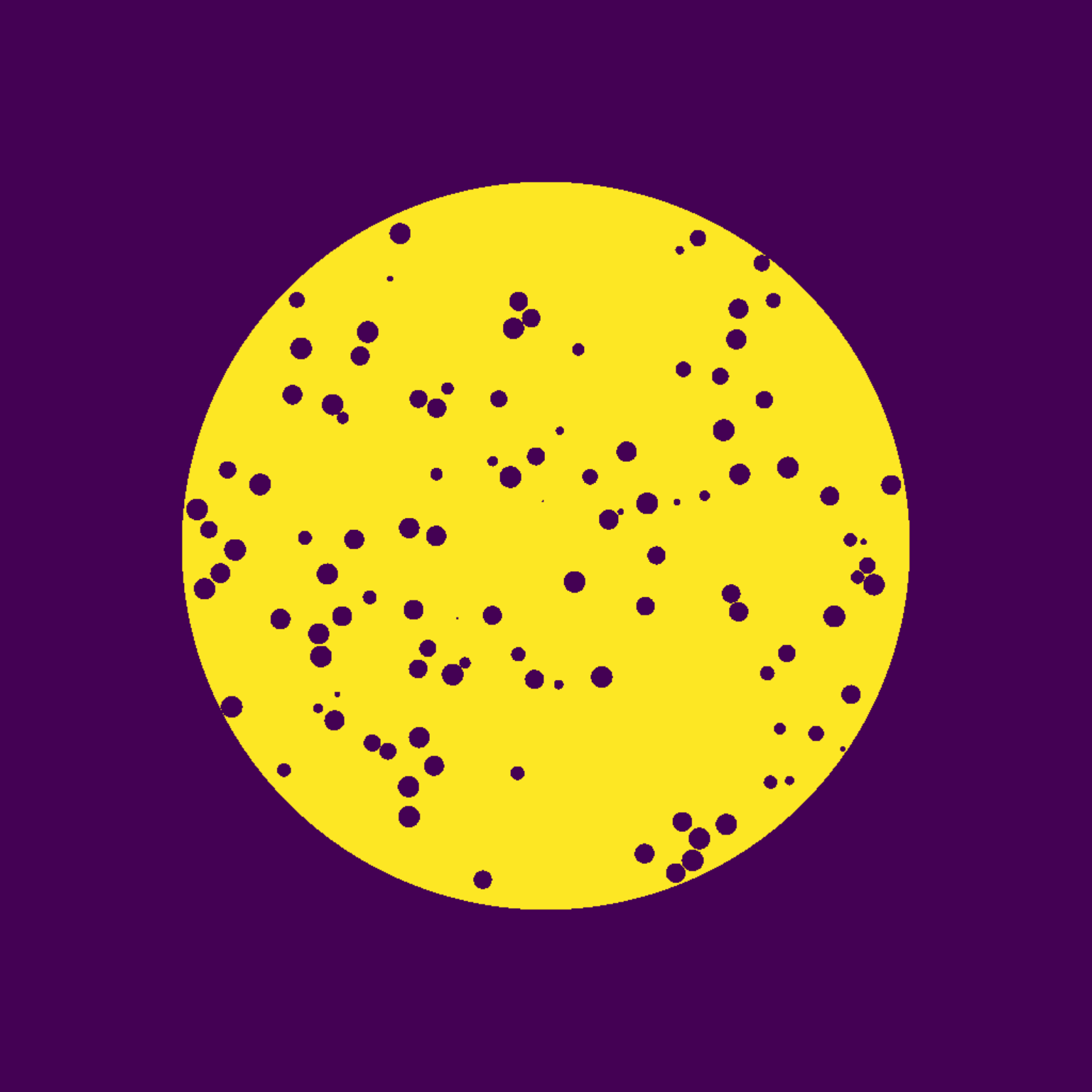}
    \includegraphics[width=0.155\textwidth]{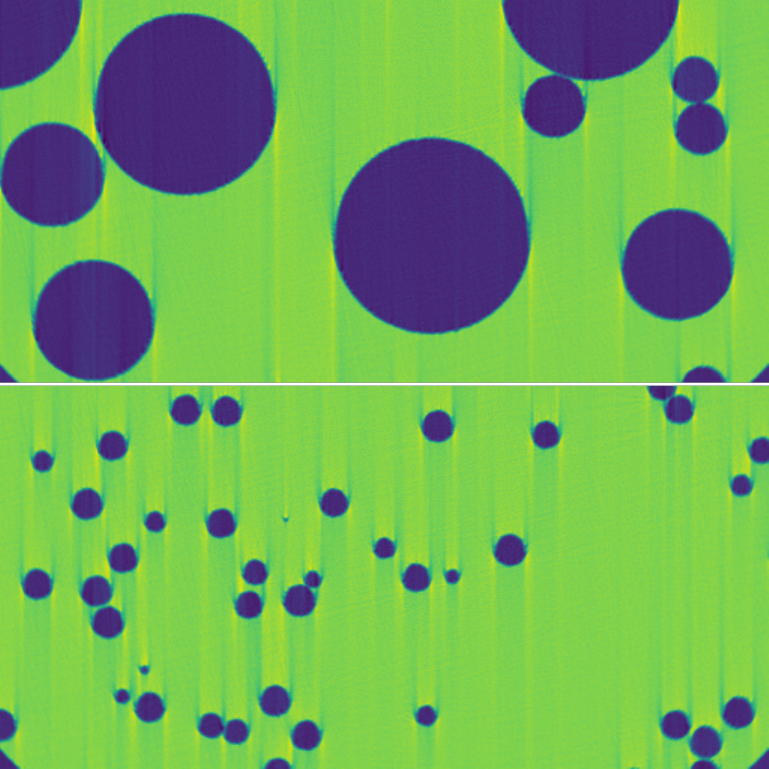}
    \caption{Phantoms, denoted \texttt{p}$_1$ (left), \texttt{p}$_2$ (middle). Right: Detail of FBP reconstructions of \texttt{p}$_1$ (top) and \texttt{p}$_2$ (bottom) with uncorrected $h=1$ (pixel) illustrating misalignment artifacts.}
    \label{fig.fan.phantoms}
\end{figure}

We simulate fan-beam projections with the \texttt{astra} toolbox \cite{astra} of two numerical phantoms based on the \texttt{foam\_ct\_phantom} library \cite{foam}, in python. The domain of the phantom is the unit disk $\{\bm x \in \R^2, \|\bm x\|_2 \leq 1\}$ sampled on a square mesh of $1024\times1024$ pixels. The projections are taken with 1024 pixels (parameterized with $s_i$) and 1024 tomographic rotations $\beta_j$ uniformly spaced with $(s_i,\beta_j) \in [-\bar s,\bar s]\times[0,2\pi)$ with $\bar s=r(r^2-1)^{-1/2}$ taken to sample the entire unit disk. Recall that $\{s_j\}$ samples the \emph{effective} pixel, \ie\ all $s_j$ are on the line crossing the origin (and the sample) and parallel to the detector, see Figure \ref{fig.fan}. The source radius is set to $r=2$. Figure \ref{fig.fan.phantoms} shows the phantoms, referred to as \texttt{p}$_1$ and \texttt{p}$_2$ with features of different sizes and the effect of 1 pixel misalignment in their FBP reconstructions. 

Sinograms are shifted in $s$ by $h^\star=10$ (effective) pixels, Yang's method \cite{yang2013} and algorithms LY, 2DR and FP are applied. The error of the recovered values $\h$ is simply expressed by $|h^\star-\h|$ (pixels). With ideal data, all four methods give almost perfect results as expected. Instead of adding random noise to data, we consider beam instabilities in both $s$ and $\beta$ by modifying the misaligned data $\g$ as in (\ref{eq.shifted}) by $\g(s,\beta) + b(s,\beta)$ for some real function $b$ modelling the beam instability. This instability was observed in our CT experiments even after flat/dark field corrections and it could be understood as the projections of the function $f\equiv 0$. Robustness to noise will be studied with real industrial CT data with different noise levels below.  

\begin{figure*}[!t]
\centering
\begin{tabular}{ll}
\begin{tikzpicture}
\begin{axis}[xlabel={$\alpha$}, ylabel={$|h^\star-\h|$ (pix.)}, legend style={at={(0,1)}, anchor=north west, font=\scriptsize}, width = .5\textwidth, height = .27\textwidth, y label style={at={(axis description cs:.05,.5)}}, xtick={0,0.002,0.004,0.006,0.008,0.01}]
\addplot coordinates {  
	(0, 0.025)
 	(0.002, 0.08)
   	(0.004, 0.18)
	(0.006, 0.25)
	(0.008, 0.32)
	(0.01, 0.49)
};
\addplot+[purple, mark=star, mark options={fill=purple}] coordinates {   
	(0, 0)
    (0.002, 0.085)
   	(0.004, 0.17)
	(0.006, 0.26)
	(0.008, 0.355)
	(0.01, 0.46)
};
\addplot+[mark=*] coordinates {  
	(0, 0.005)
    (0.002, 0.025)
   	(0.004, 0.055)
 	(0.006, 0.085)
	(0.008, 0.115)
	(0.01, 0.145)
};
\addplot coordinates {  
	(0, 0)
    (0.002, 0.025)
   	(0.004, 0.055)
	(0.006, 0.08)
	(0.008, 0.11)
	(0.01, 0.135)
};
    \legend{Yang, LY, FP, 2DR}
	\end{axis}
\end{tikzpicture}  &
\begin{tikzpicture}
\begin{axis}[xlabel={$\alpha$}, ylabel={$|h^\star-\h|$ (pix.)}, legend style={at={(0,1)}, anchor=north west, font=\scriptsize}, width = .5\textwidth, height = .27\textwidth, y label style={at={(axis description cs:.05,.5)}}, xtick={0,0.002,0.004,0.006,0.008,0.01}]
\addplot coordinates {  
	(0, 0.01)
 	(0.002, 0.125)
   	(0.004, 0.23)
	(0.006, 0.3)
	(0.008, 0.45)
	(0.01, 0.555)
};
\addplot+[purple, mark=star, mark options={fill=purple}] coordinates {   
	(0, 0)
    (0.002, 0.1)
   	(0.004, 0.21)
	(0.006, 0.32)
	(0.008, 0.445)
	(0.01, 0.58)
};
\addplot+[mark=*] coordinates {  
	(0, 0.005)
    (0.002, 0.06)
   	(0.004, 0.115)
 	(0.006, 0.17)
	(0.008, 0.22)
	(0.01, 0.275)
};
\addplot coordinates {  
	(0, 0)
    (0.002, 0.045)
   	(0.004, 0.09)
	(0.006, 0.135)
	(0.008, 0.18)
	(0.01, 0.23)
};
    \legend{Yang, LY, FP, 2DR}
	\end{axis}
\end{tikzpicture} 
\end{tabular}
\caption{Shift estimation results for phantoms \texttt{p}$_1$ (left) and \texttt{p}$_2$ (right), showing better performances of 2DR and FP when the beam instability increases.}
\label{fig.fan.simulated}
\end{figure*}
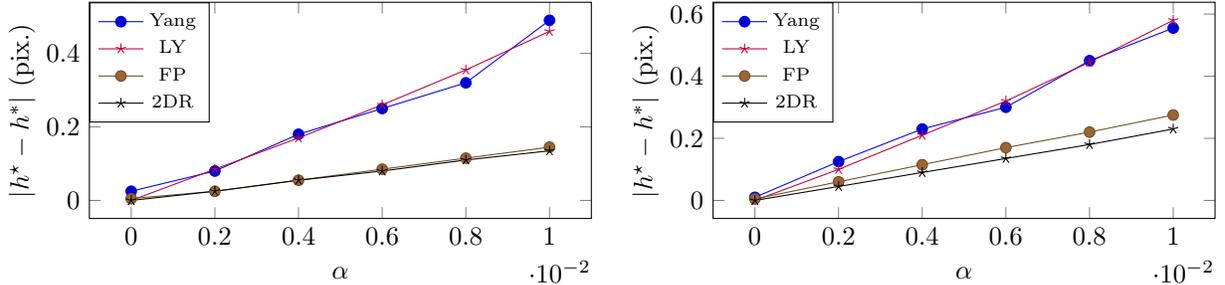

For simplicity, we consider $b(s,\beta) = \alpha ( \sin(\pi s / (2\bar s)) + \cos(\replaced{\beta}{\theta} / 2) +2 )$, then the instabilities are positive, monotonic with sinusoidal behaviour of half period \added{in both variables}. $\alpha>0$ controls the amount of instability. Figure \ref{fig.fan.simulated} shows the results on both phantoms \texttt{p}$_1$ and \texttt{p}$_2$ for increasing values of $\alpha$. With $\alpha=0$, exact values are recovered with algorithms LY, 2DR and FP, and with negligible error for Yang's method (0.02 of a pixel). By increasing $\alpha$ the two proposed algorithms (2DR and FP) obtain significantly better results. The reason is what we mentioned in the presentation of the methods, Yang and LY \emph{average out} the data by the integral operation and then the image registration relies mainly on the boundary of the projections while 2DR and FP, by avoiding any averaging, incorporate all internal features in the registration. The improvement will be more obvious in the next section with \replaced{experimental}{real} data. Python codes are available at \texttt{github.com/patoguerrero/alignCT} to reproduce the case $\alpha=0$ and are easily adaptable to any experimental fan-beam data. 
\begin{table}
\centering \small 
\caption{CT settings for the acquired industrial data.}
\begin{tabular}{lllll}
\addlinespace[0pt]\cmidrule{1-5}\addlinespace[1pt]
 & voltage & current & views & time / view \\    
\texttt{low noise}  & 130 kV. & 140 µA. & 3142 & 1415 µs. \\    
\texttt{high noise}  & 110 kV. & 127 µA. & 720 & 500 µs. \\    
\addlinespace[1pt]\cmidrule{1-5}
\end{tabular}
\label{tab.settings}
\end{table}
\begin{table*}
\centering \small 
\caption{MSE and $\h$ results for 2D industrial data with two different noise levels.}
\begin{tabular}{lllllllllll}
\addlinespace[0pt]\cmidrule{1-11}\addlinespace[1pt]
& \multicolumn{5}{l}{\texttt{low noise}} & \multicolumn{5}{l}{\texttt{high noise}} \\ 
\addlinespace[2pt]\cmidrule{2-11}\addlinespace[2pt]
& Yang & LY & FP & FP$_{10}$ & 2DR &  Yang & LY & FP & FP$_{10}$ & 2DR \\    
$\h$ & 7.05 & 7.06 & 5.08 & 4.45 & 4.45 & 3.54 & 3.5 & 2.02 & 1.89 & 1.98 \\    
MSE  & 0.5 & 0.5 & \textbf{0.419} & \textbf{0.41} & \textbf{0.41} & 0.253 & 0.254 & \textbf{0.196} & \textbf{0.2} & \textbf{0.2} \\    
\addlinespace[1pt]\cmidrule{1-11}
\end{tabular}
\label{tab.fan}
\end{table*}
\begin{figure*}[!t]
    \centering
    \includegraphics[height=0.22\textwidth]{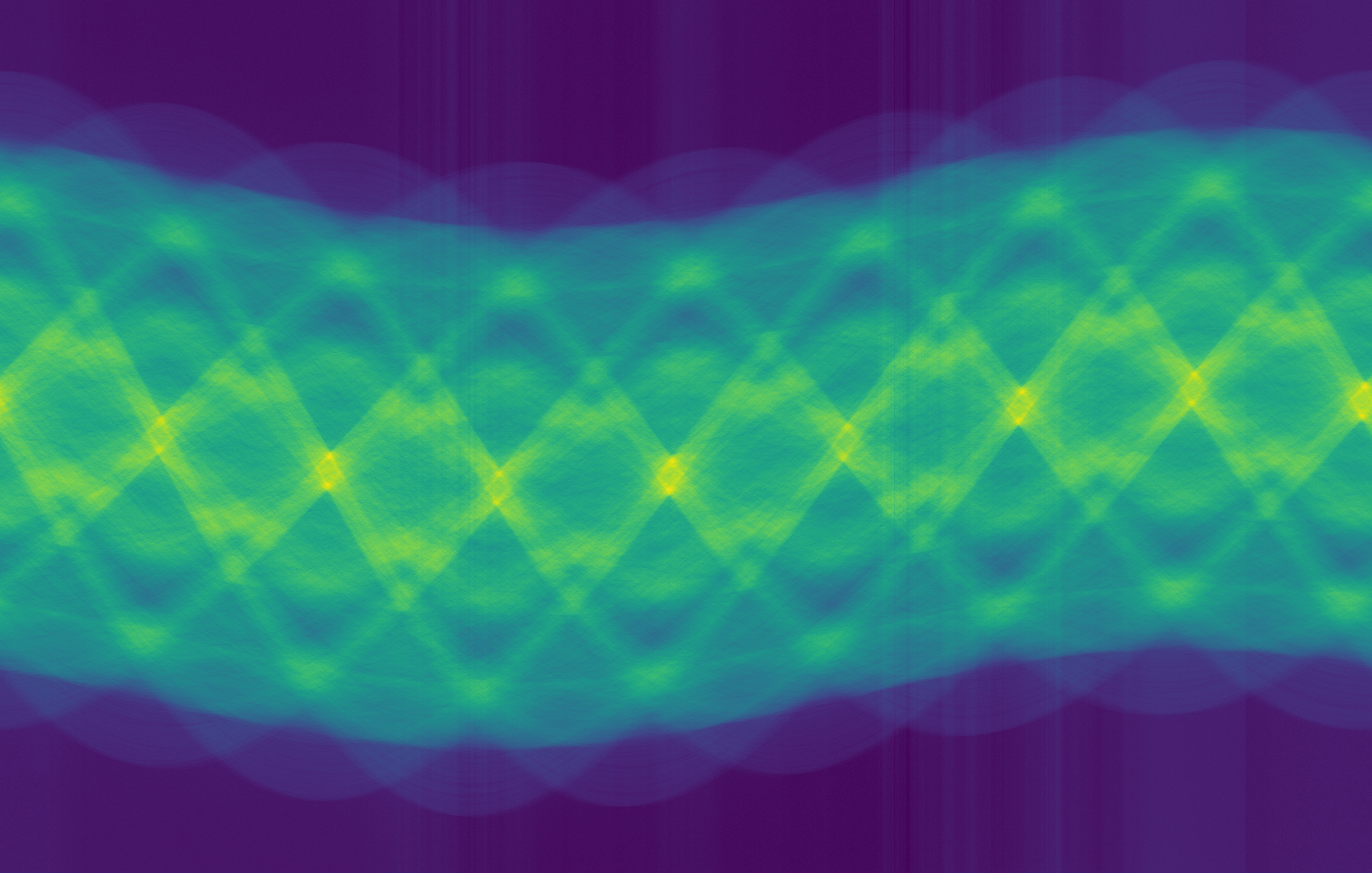}
    \put(-0.36\textwidth,0.1\textwidth){\small $s$} \put(-0.2\textwidth,-0.025\textwidth){\small $\beta \in [0,2\pi)$} \ 
    \includegraphics[height=0.22\textwidth]{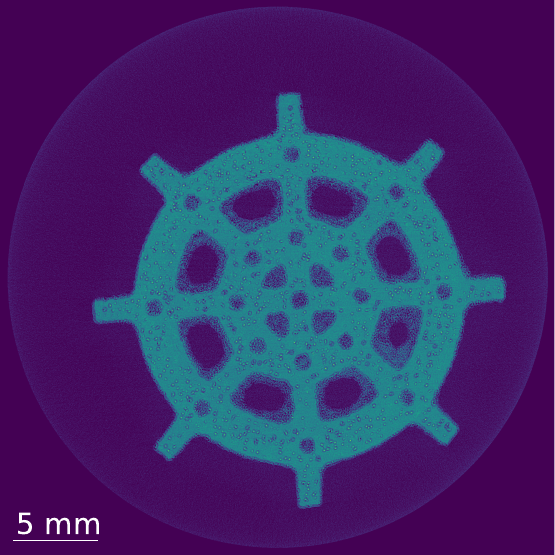}
    \includegraphics[height=0.22\textwidth]{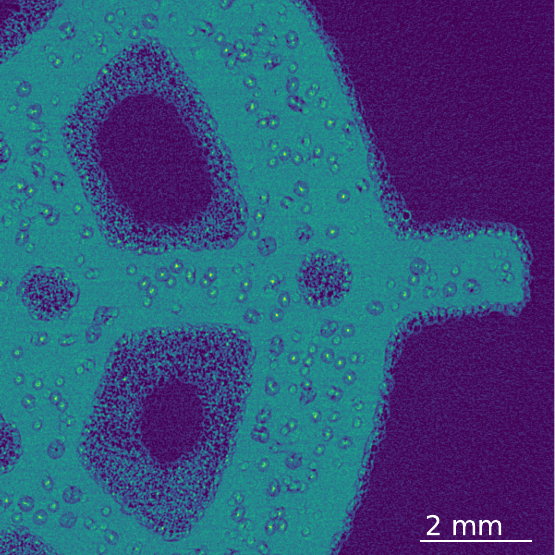}
    \caption{\texttt{low noise} fan-beam sinogram of the manufactured object (left), followed by its FBP reconstruction without any detector shift estimation and its detail where double edge artifacts due to misalignments are visible.}
    \label{fig.fan.timonHQ0}
\end{figure*}
\begin{figure*}[!t]
    \centering
    \rotatebox{90}{\small LY\phantom{$_{10}$}}
    \includegraphics[width=0.22\textwidth]{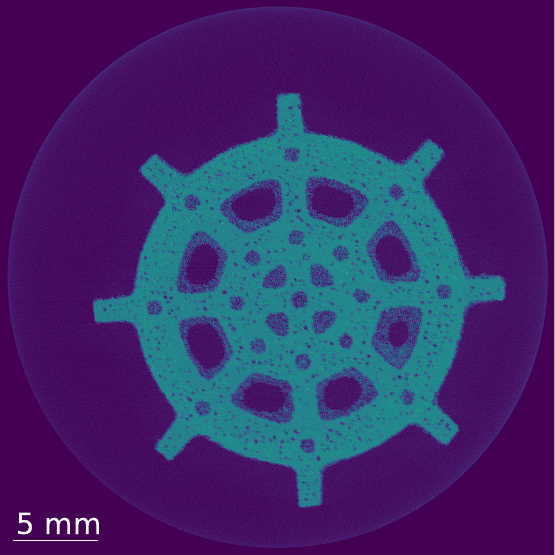}
    \includegraphics[width=0.22\textwidth]{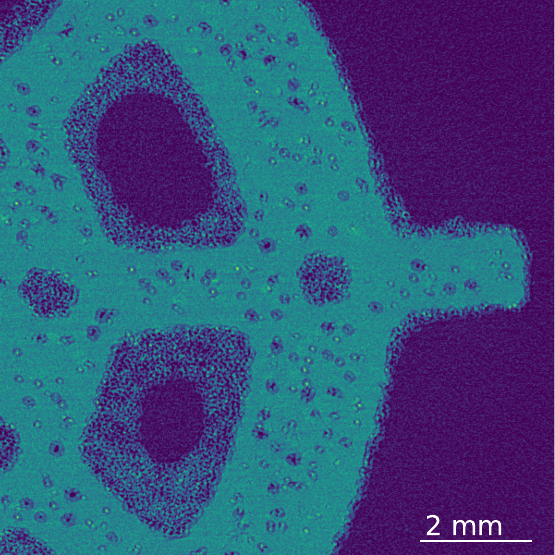}
    \includegraphics[width=0.22\textwidth]{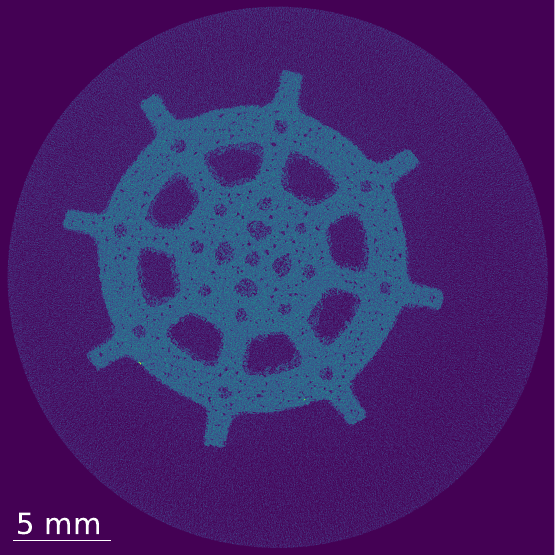}
    \includegraphics[width=0.22\textwidth]{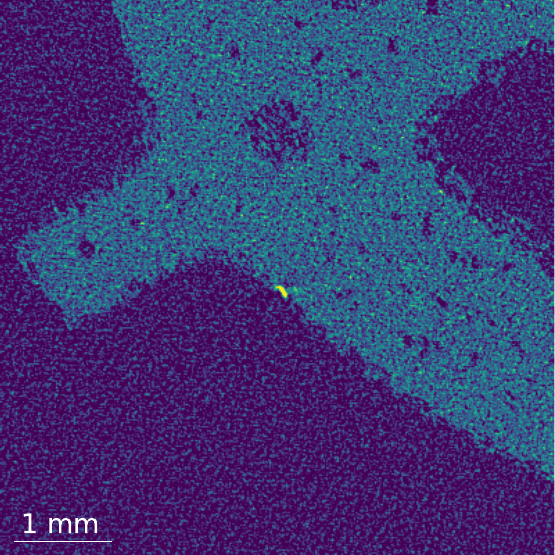} \\
    \rotatebox{90}{\small FP\phantom{$_{10}$}}
    \includegraphics[width=0.22\textwidth]{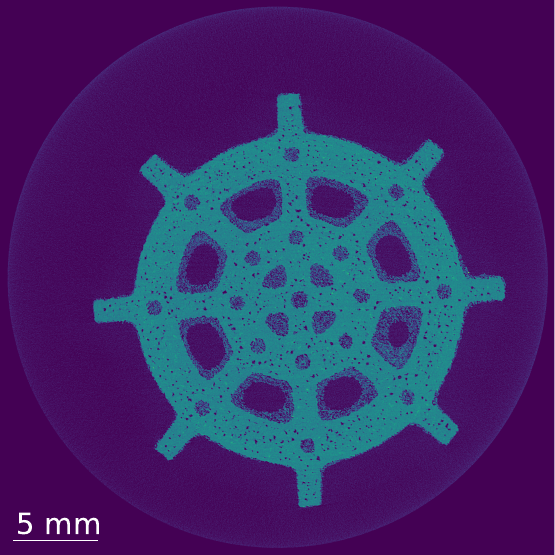}
    \includegraphics[width=0.22\textwidth]{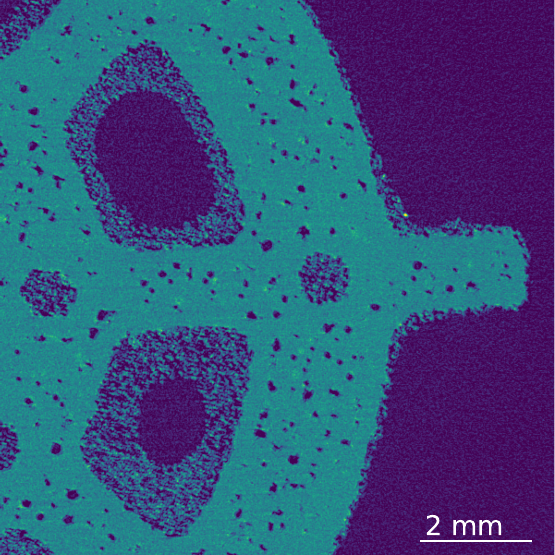}
    \includegraphics[width=0.22\textwidth]{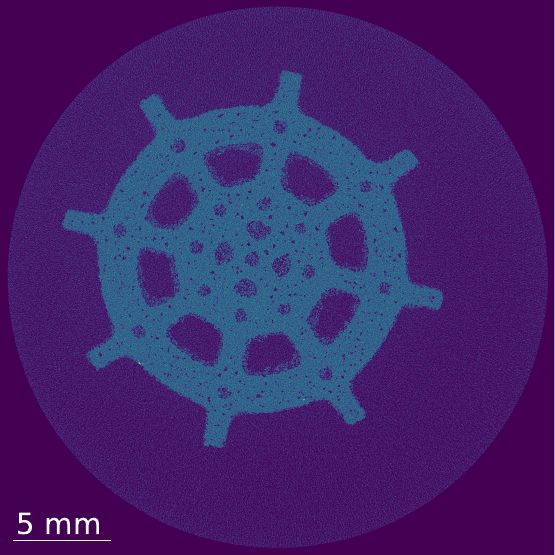}
    \includegraphics[width=0.22\textwidth]{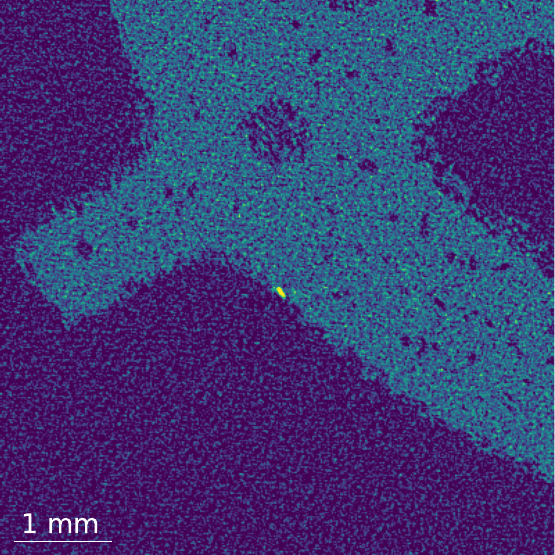} \\
    \rotatebox{90}{\small FP$_{10}$}
    \includegraphics[width=0.22\textwidth]{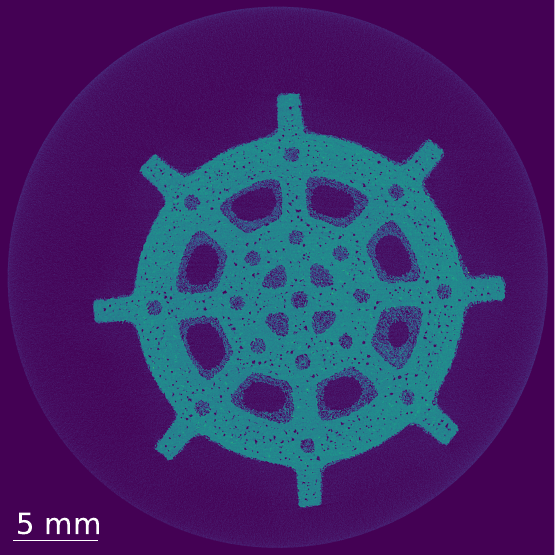} 
    \includegraphics[width=0.22\textwidth]{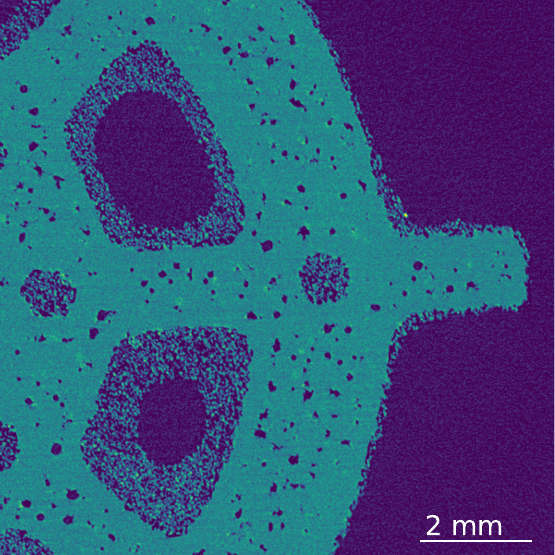} 
    \includegraphics[width=0.22\textwidth]{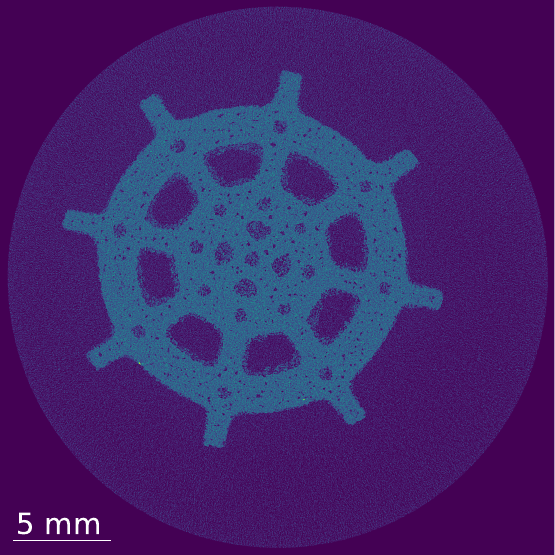} 
    \includegraphics[width=0.22\textwidth]{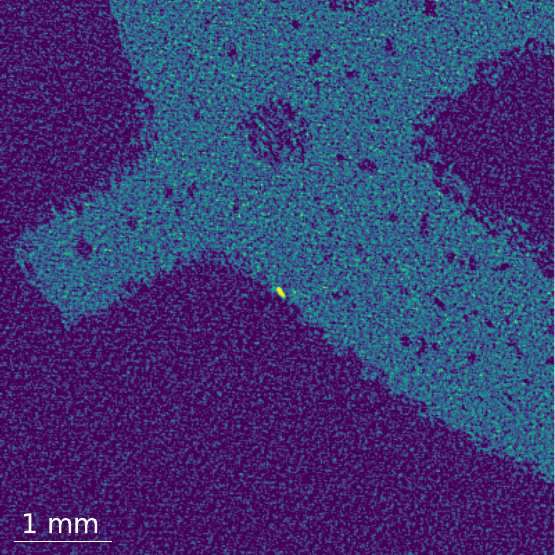}\\
    \texttt{low noise} \hspace{0.32\textwidth} \texttt{high noise}
    \caption{FBP reconstructions of industrial experiments \texttt{low noise} and \texttt{high noise} with their details. Yang's method and 2DR are not shown because they are similar to LY and FP$_{10}$ respectively. Strong artifacts visible with LY are corrected with FP and slightly better with FP$_{10}$.}
    \label{fig.fan.timon}
\end{figure*}

\subsubsection{Industrial CT \replaced{experiment}{data}} A Nikon XT H 225 ST CT scanner that operates in cone-beam geometry was used, fan-beam sinograms are simply taken as the middle row of the projections. As no ground-truth data are available, results will be quantified based in relationship (\ref{eq.shiftedinq}) as the methods were derived from this necessary condition. Namely, we use the normalized mean squared error (MSE) computed by $\|\g - \hat g \|_2^2 / \|\g\|_2^2$ between the data $\g$ and its \emph{symmetric} sinogram $\hat g\colon (s,\beta)\mapsto \g (-s+2\h, \beta+\pi+2\arctan\left((s-\h)/r\right))$, with $\|\cdot\|_2$ the Hilbert–Schmidt norm, or the Frobenius norm in the discrete setting.

A designed and additive manufactured object printed with laser sintering of Polyamide 12 was CT scanned with two different scanning settings, listed in Table \ref{tab.settings}. These configurations were established to produce low and high noise level experiments. We refer to \texttt{low noise} the former experiment and to \texttt{high noise} the latter. Note that the total scanning time for \texttt{low noise} was 75 min. while for \texttt{high noise} it was 6 min. Polymer laser sintered material has an internal pore morphology with pores (observed in the reconstruction images) in the range of 30 to 300 microns.

Figure \ref{fig.fan.timonHQ0} shows the central fan-beam sinogram of \texttt{low noise} and their FBP reconstruction without any detector shift correction, clearly no analysis can be done with such reconstructions due to the observed double edge artifacts. Then the \replaced{5}{4} algorithms (Yang, LY, 2DR, FP, FP$_{K}$) are used to estimate $\h$ for both high/low noise level experiments \deleted{and an additional FP$_{10}$ method explained below}. Table \ref{tab.fan} presents the estimated values and their corresponding MSE. As expected, we observe that FP and 2DR significantly outperform  Yang's and LY. 2DR is the best method overall as expected as well. Note that the performance of Yang's and LY is almost identical suggesting than Yang's idea of averaging data cannot be improved by higher order interpolations. The performance of 2DR and FP is however, close enough given their difference in computational cost. With \texttt{low noise}, FP and 2DR differ in near half a pixel. This is because FP performs 1D signal registration by fixing $\beta=0$ in the signal $\Lambda$ and computes \added{iteratively} its \emph{approximately symmetrical} signal $\Pi$ (see Algorithm \ref{alg.fp}). \added{Nevertheless, FP$_{10}$ and 2DR give equivalent results, recalling than $M=5$ in all experiments after verifying that convergence is reached, thus FP$_{10}$ still has lower cost than 2DR.} 

\deleted{We can improve the robustness of FP by simply taking the median of \eg, $K=10$ calls of FP with different values of $\beta$ uniformly distributed in the angular range, then an eventual defective projection at $\beta=0$ will have no effect. This modified algorithm, denoted by FP$_{K}$ still has lower cost than 2DR as long as $KM\ll N_\beta$. Recall than here $M=5$ in all experiments after verifying that convergence is reached. Note that FP$_{10}$ and 2DR give now equivalent results.} 

\begin{figure}[!t]
\centering
\rotatebox{90}{\small\texttt{low noise}\phantom{g}}
\subfloat[LY]{\includegraphics[width=0.1\textwidth]{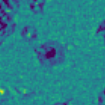}}\
\subfloat[FP]{\includegraphics[width=0.1\textwidth]{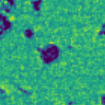}}\
\subfloat[FP$_{10}$]{\includegraphics[width=0.1\textwidth]{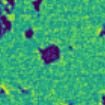}}  \\ 
\rotatebox{90}{\small\texttt{high noise}}
\subfloat[$\h = 0$]{\includegraphics[width=0.1\textwidth]{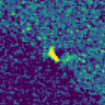}}\
\subfloat[LY]{\includegraphics[width=0.1\textwidth]{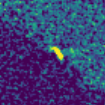}}\
\subfloat[FP]{\includegraphics[width=0.1\textwidth]{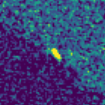}} 
\caption{Porous detail of reconstructions of \texttt{low noise} obtained with (a) LY, (b) FP and (c) FP$_{10}$; contamination detail of reconstructions of \texttt{high noise} obtained with (d) $\h = 0$ (no correction), (e) LY and (f) FP.}
\label{fig.fan.porous}
\end{figure}

Figure \ref{fig.fan.timon} shows FBP reconstruction images with LY, FP and FP$_{10}$ shift estimation. Yang's and 2DR images are not presented because we could not differentiate them visually from LY and FP$_{10}$ respectively, as expected given the values in Table \ref{tab.fan}. Strong artifacts are still present with LY, while small artifacts are identifiable at the porous level with FP but not with FP$_{10}$, as exposed in Figure \ref{fig.fan.porous}(a,b,c). 

Finally, we can compare visually the performance of the methods with highly noisy data by observing the small highly attenuating region (in yellow) that represents a metallic contamination in the printed sample. Yang and LY over-estimate $\h \approx 3.5$ pixels while FP/FP$_{10}$/2DR obtain $\h\approx2$ pixels. The particle with $\h = 0$ (no correction), LY and FP is exposed in Figure \ref{fig.fan.porous}(d,e,f) where we observe the over-estimating effect of LY as an ellipse-like particle seems to be deformed in opposite directions for $\h = 0$ and LY with respect to FP. 

To conclude this section, we observed that Yang's method can be theoretically improved by LY that incorporates an interpolation degree other than plain nearest neighbour. However, the improvement is not enough as the main drawback of averaging the data is still present. 2DR and FP appear to solve this issue. To our knowledge, LY, 2DR and FP are presented for the first time.


\subsection{\added{Cone-beam data}}
\subsubsection{\added{Simulated experiment}} 
\begin{figure}[!t]
    \centering
    \includegraphics[width=0.155\textwidth]{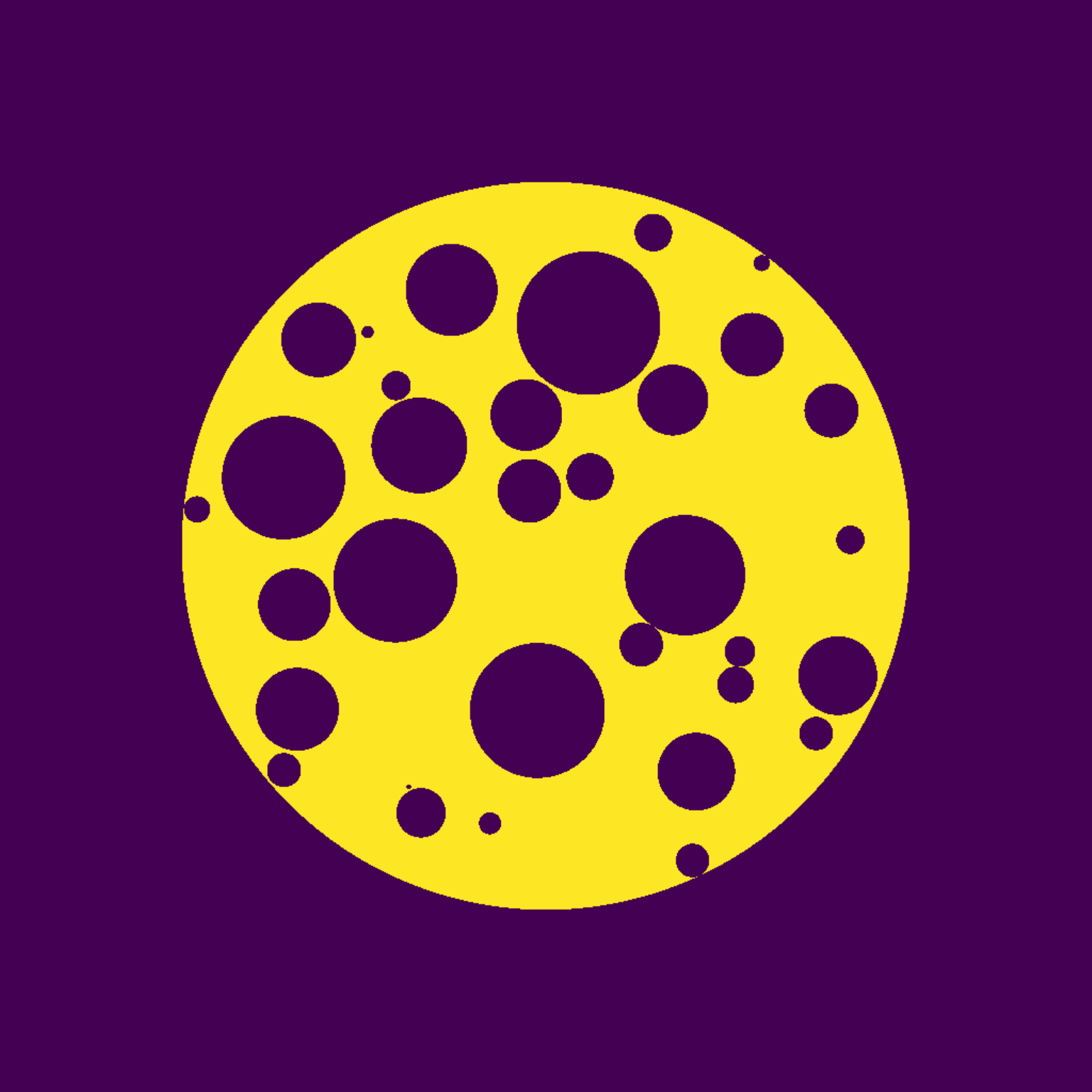}
    \includegraphics[width=0.155\textwidth]{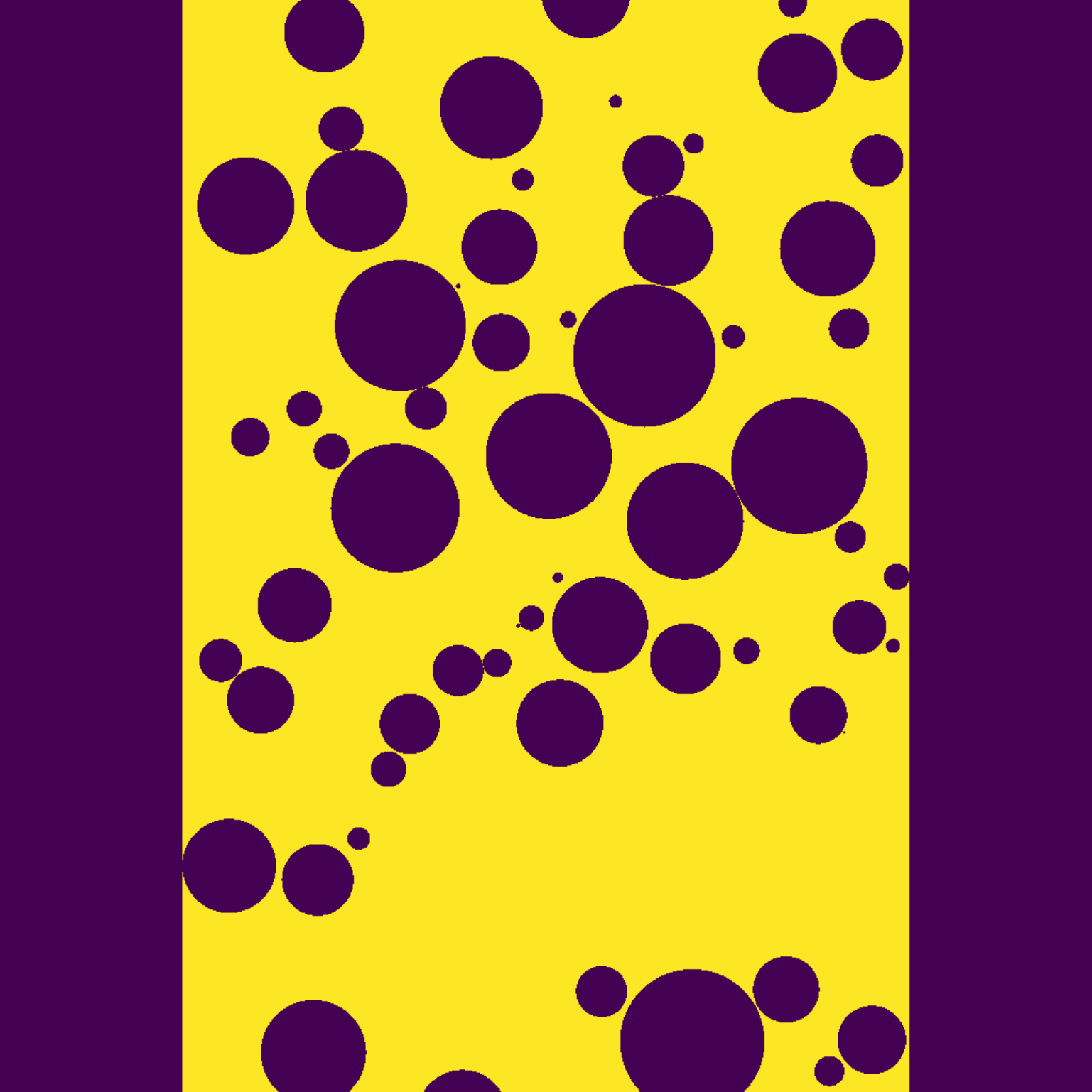}
    \includegraphics[width=0.155\textwidth]{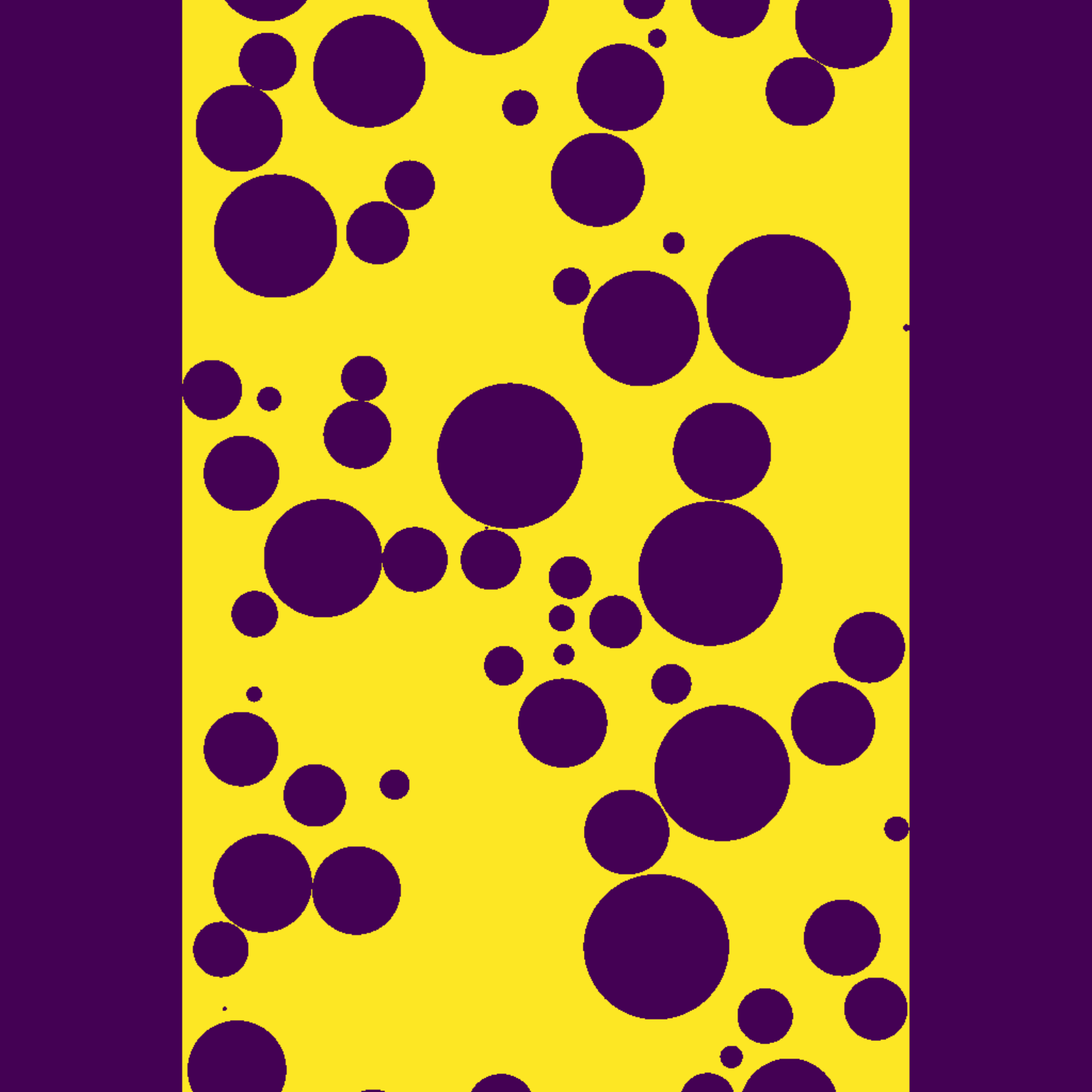}
    \caption{Orthogonal central slices of the simulated phantom, $y=0$, $x=0$, $z=0$ respectively.}
    \label{fig.cone.simu}
\end{figure}

The following simulation experiment can be reproduced from \texttt{github.com/patoguerrero/alignCT}. The \texttt{foam\_ct\_phantom} python library is used to generate a 3D cylindrical phantom composed of void spheres with slices similar to phantom \texttt{p}$_1$ of Section \ref{sec.2d.results}. The orthogonal central slices of the phantom are displayed in Figure \ref{fig.cone.simu}. The phantom (and projection data) has $1024\times1024\times1024$ voxels with equivalent sampling and configurations of Section \ref{sec.2d.results}. Cone-beam projections were obtained with the \texttt{astra} toolbox and misaligned with $h^\star = 10$ (effective) pixels and $\eta^\star = 1$ degree. The VP approach was executed with both FP$_{10}$ and 2DR to solve the inner problem $\textsf{h}(\gb, \eta_k)$, referred to respectively as VP-FP$_{10}$ and VP-2DR. Both methods obtained satisfactory results in 2 gradient descent iterations. Namely, $(\h,\et) = (9.98, 1.0192)$ with VP-FP$_{10}$ and $(\h,\et) = (10, 1.0196)$ with VP-2DR.  

\subsubsection{Industrial CT \replaced{experiment}{data}} 
\begin{figure}
    \centering
    \includegraphics[width=0.17\textwidth]{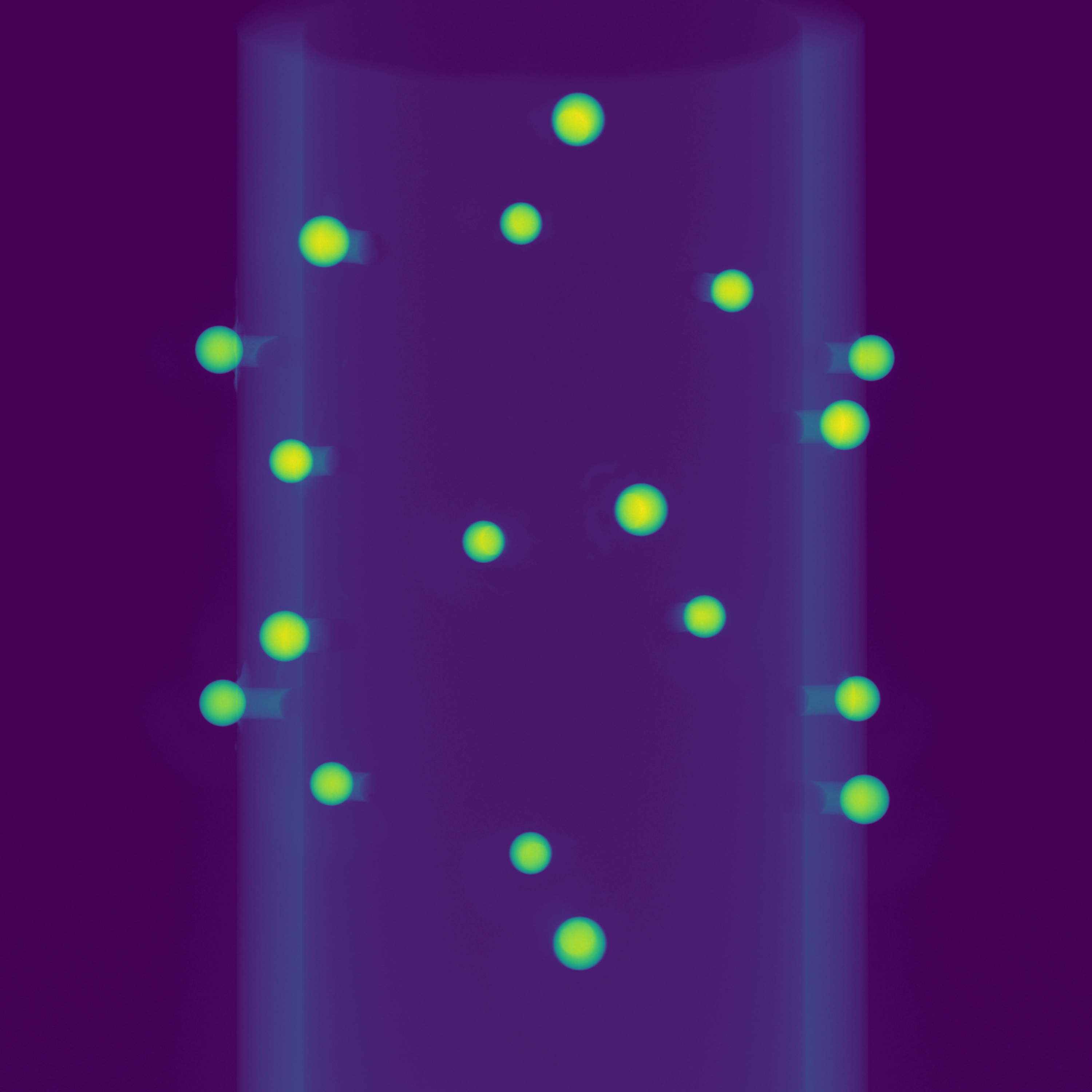}
    \includegraphics[width=0.17\textwidth]{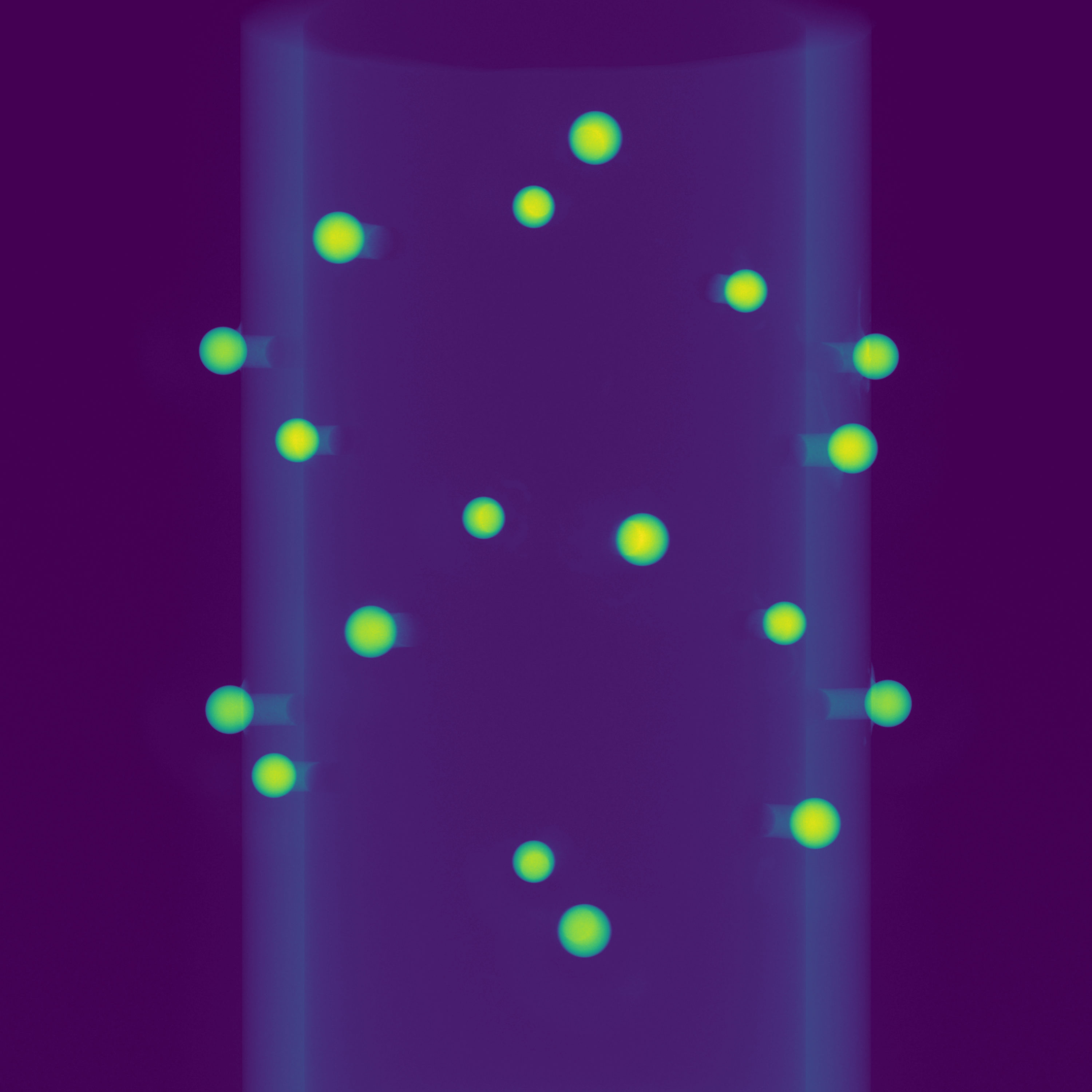}
    \caption{Projections of the \texttt{calibration} object under two views separated by $\pi$ radians.}
    \label{fig.cone.gepettodata}
\end{figure}

\begin{figure*}[!t]
    \centering
    \rotatebox{90}{\azul{\texttt{calibration\phantom{g}}(raw)}}
    \includegraphics[width=0.22\textwidth]{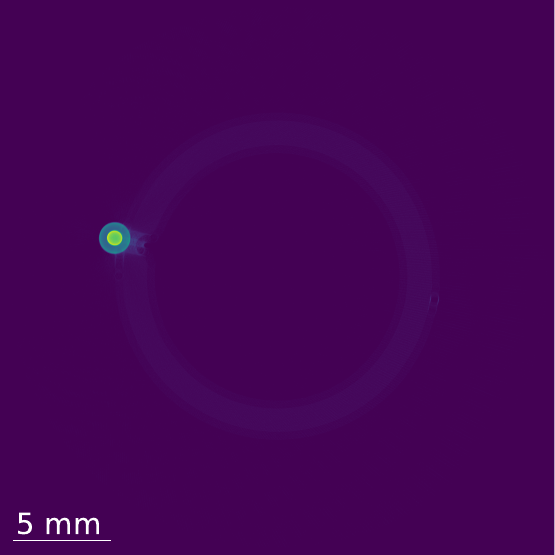}
    \includegraphics[width=0.22\textwidth]{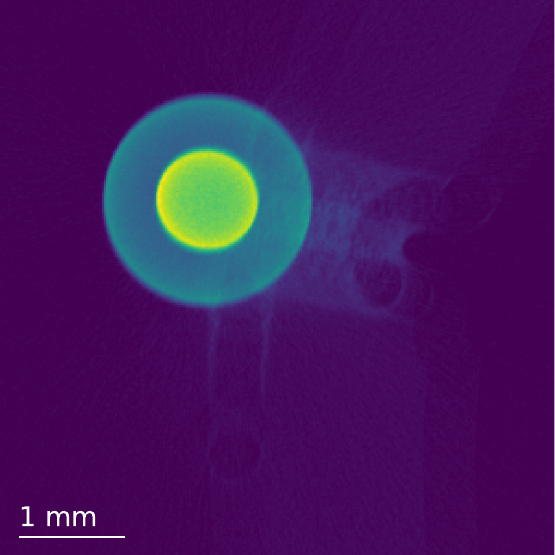}
    \includegraphics[width=0.22\textwidth]{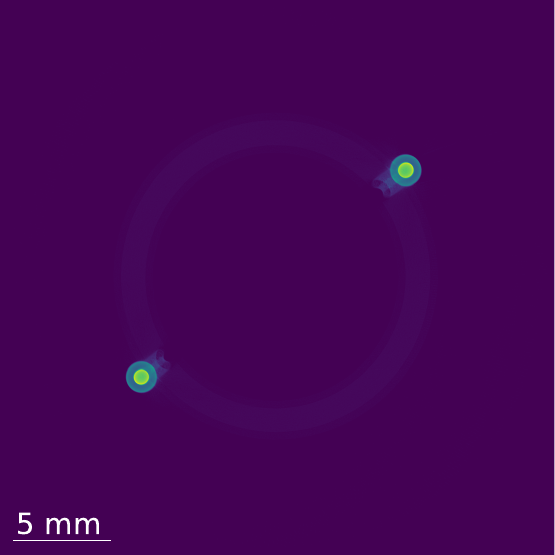}
    \includegraphics[width=0.22\textwidth]{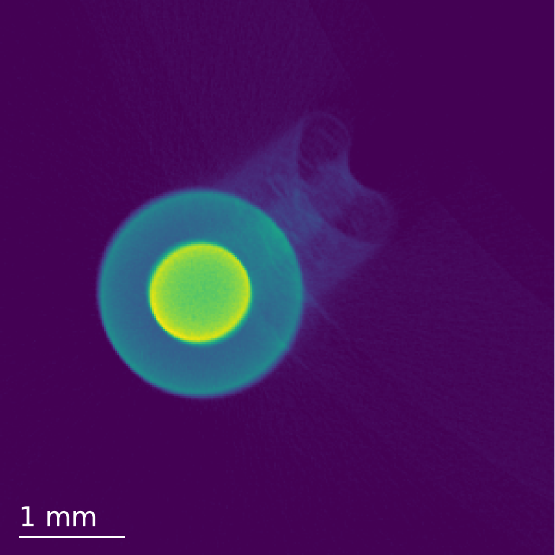} \\
    \rotatebox{90}{\texttt{calibration\phantom{g}}}
    \includegraphics[width=0.22\textwidth]{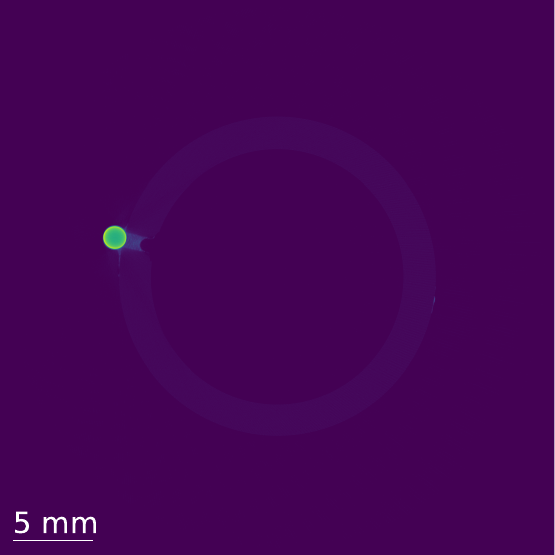}
    \includegraphics[width=0.22\textwidth]{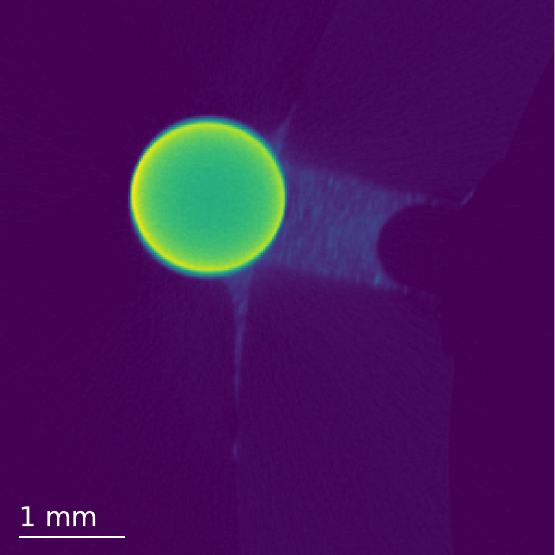}
    \includegraphics[width=0.22\textwidth]{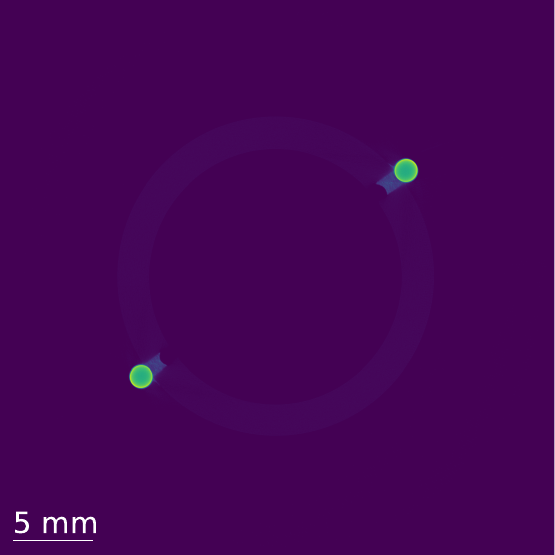}
    \includegraphics[width=0.22\textwidth]{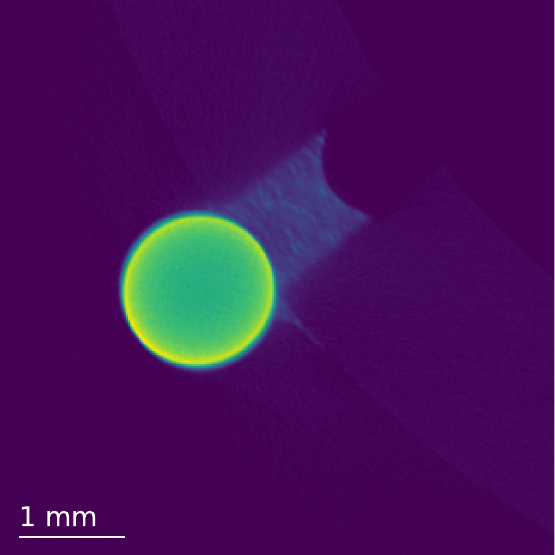} \\
    \rotatebox{90}{\texttt{low noise\phantom{g}}}
    \includegraphics[width=0.22\textwidth]{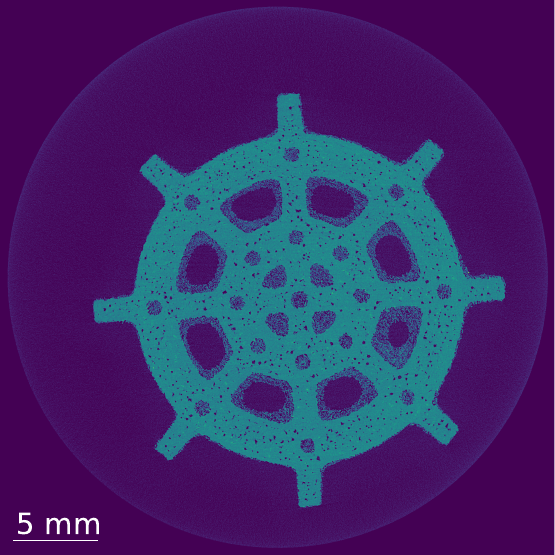}
    \includegraphics[width=0.22\textwidth]{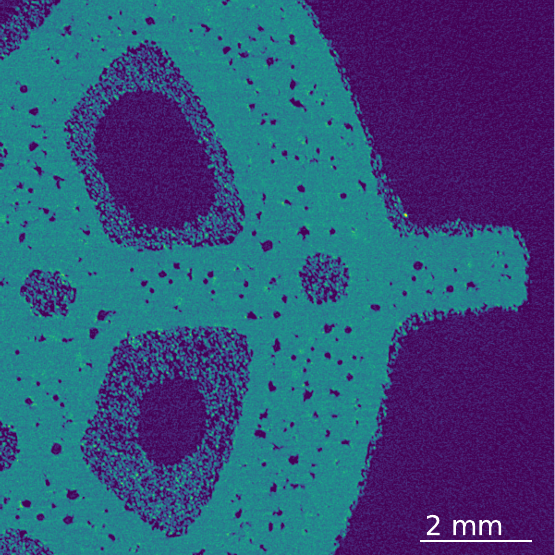}
    \includegraphics[width=0.22\textwidth]{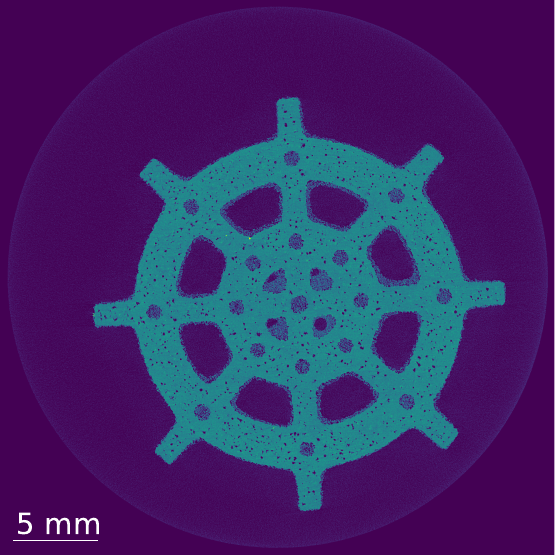}
    \includegraphics[width=0.22\textwidth]{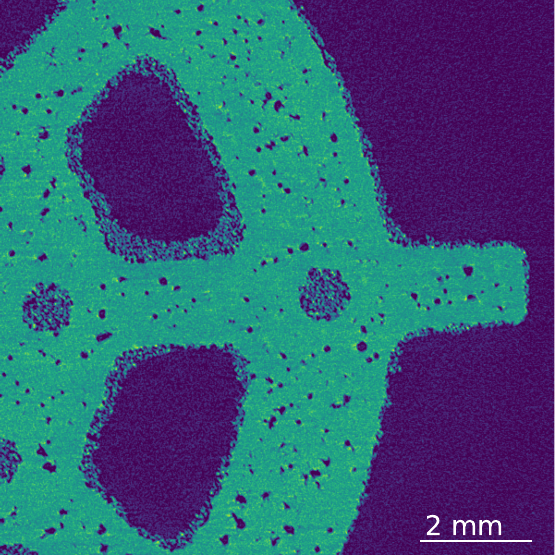} \\
    \rotatebox{90}{\texttt{high noise}}
    \includegraphics[width=0.22\textwidth]{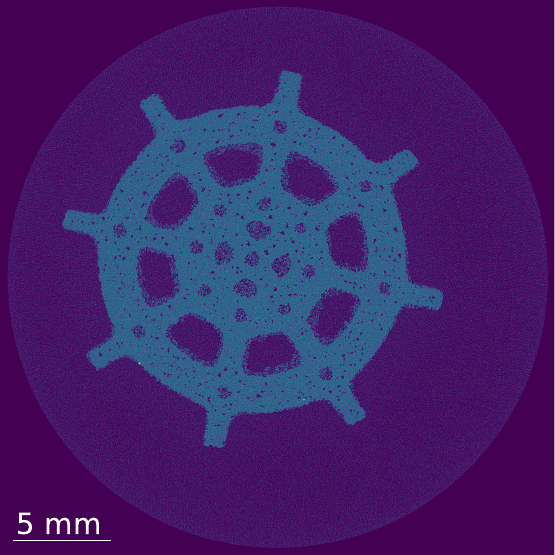}
    \includegraphics[width=0.22\textwidth]{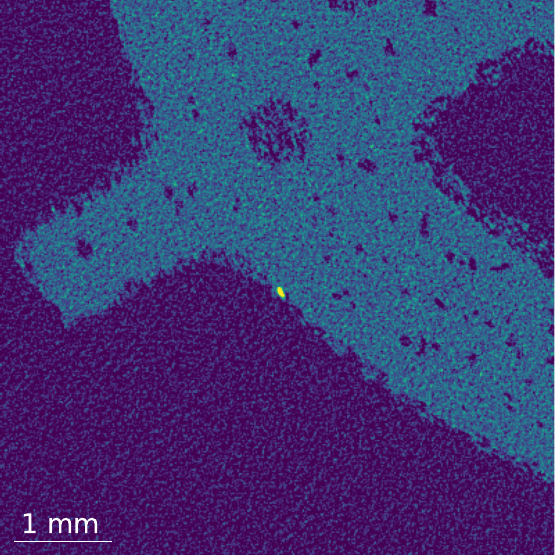}
    \includegraphics[width=0.22\textwidth]{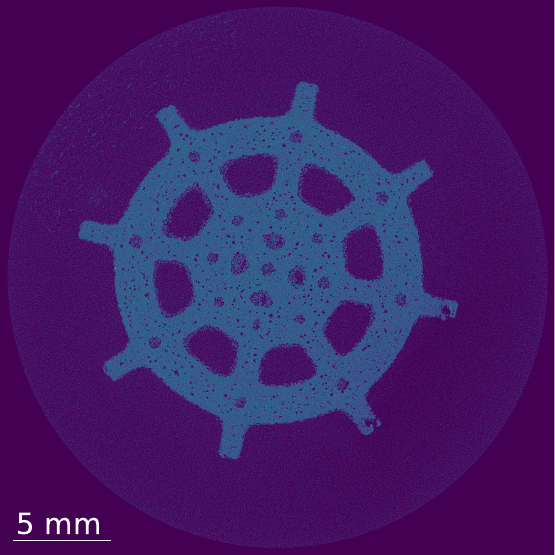}
    \includegraphics[width=0.22\textwidth]{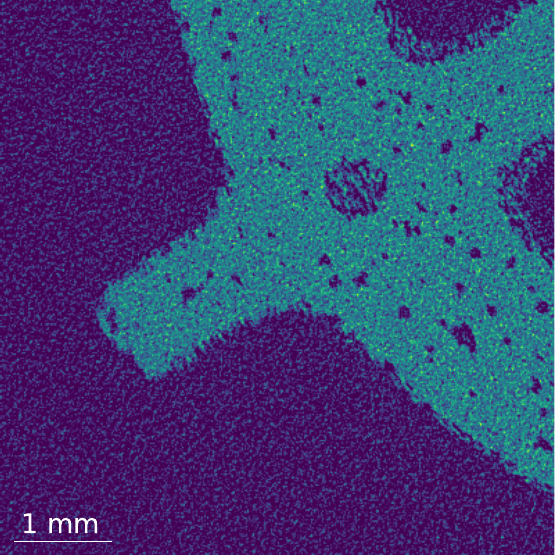}\\
    $y=0$ \hspace{0.34\textwidth} $y=512$ 
    \caption{FDK reconstructions of the \texttt{calibration} object and of the additive manufactured object with low and high noise levels with their details. The alignment is done with VP-2DR, no visual difference is observed with respect to VP-FP$_{10}$ reconstructions (not shown). \added{The top row shows \texttt{calibration} with no alignment compensation exhibiting double edge artifacts.} Note \added{also} the presence of the sharp low-attenuating region attached to the metallic spheres in \texttt{calibration} \added{after alignment}, it corresponds to the used glue to fix the spheres.}
    \label{fig.cone.timon}
\end{figure*}

The same Nikon CT scanner of Section \ref{sec.2d.results} is used in the following experiments. First, the object used to estimate the CT geometry in \cite{massi} composed of steel spheres attached to a carbon fiber cylinder is measured and referred to \texttt{calibration}. We applied such off-line method also described in \cite{massi} for comparison purposes. Note that the values obtained with that method are supposed to be used in a subsequent scan and then subject to non-reproducible errors, while here we used them in the same experiment. Therefore very good results are expected. 720 projections were measured of this object, Figure \ref{fig.cone.gepettodata} shows two acquired images with views separated by $\pi$ radians after applying the linearization given by the Beer-Lambert law \cite{kak}. Table \ref{tab.cone} shows the results obtained with the reference-object method of \cite{massi} (where spheres dimensions and distances are known a-prori) and the proposed approaches VP-FP$_{10}$ and VP-2DR. The high accuracy of the later is clear with respect to the references values of \texttt{calibration}. \azul{Here, we observe that the obtained detector in-plane rotation angle $\et$ is very close to 0 (0.055 deg.) which is also the initial value $\eta_0$ used in the gradient descent iteration of Algorithm \ref{alg.vp}. Thus, we decided to test the case when $\eta$ is not very close to the ideal geometry and set the initial value as $\eta_0=-1$ deg. The obtained results are reported in Table \ref{tab.cone} and still show good correspondence with the reference values. Interestingly, with VP-FP$_{10}$ the $\et$ output is even closer to the reference value compared when $\eta_0=0$. This is because the algorithm performed 3 gradient descent iterations with 5 total backtracking iterations for the step sizes until reaching convergence, while with $\eta_0=0$ only one gradient descent iteration with 3 backtracking iterations was needed.} Finally the cone-beam scans of the additive manufactured object with low and high noise levels presented in Section \ref{sec.2d.results} are used here and the results are at Table \ref{tab.cone} as well. Note that these scans were performed after a different calibration of the machine with more accuracy. Therefore the estimated values of $\h$ are lower than 1 mm., and $\et$ is estimated as 0 in both. Convergence of the algorithms were always reached before 5 gradient descent iterations. Two 2D slices (middle height $y=0$ and $1/4$ height $y=512$ pix.) of these three mentioned experiments are showed in Figure \ref{fig.cone.timon} only with VP-2DR alignment as no visual difference is observed related to VP-FP$_{10}$. \azul{In the same figure, \texttt{calibration} without any alignment is also exposed, illustrating the severity and kind of artifacts the methods are addressing}. The obtained reconstruction seem not to present any misalignment (double edge, distortion) artifacts. The showed FDK reconstructions were performed with the \texttt{astra} toolbox with a custom defined geometry incorporating $(\h,\et)$ in the backprojection step.   

\begin{table*}
\centering \small 
\caption{Results for 3D industrial data. \azul{Columns after VP-FP$_{10}$ and VP-2DR for \texttt{calibration} report results with an initial value $\eta_0$ in Algorithm \ref{alg.vp} further away from the reference value (in deg.).}}
\begin{tabular}{llllllllll}
\addlinespace[4pt]\cmidrule{1-10}\addlinespace[1pt]
& \multicolumn{5}{l}{\texttt{calibration}} & \multicolumn{2}{l}{\texttt{low noise}} & \multicolumn{2}{l}{\texttt{high noise}} \\ 
\addlinespace[2pt]\cmidrule{2-10}\addlinespace[2pt]
& reference & VP-FP$_{10}$ & \azul{$\eta_0 = -1$}  & VP-2DR & \azul{$\eta_0 = -1$} & VP-FP$_{10}$ & VP-2DR & VP-FP$_{10}$ & VP-2DR \\    
$\h$ (mm.) & 2.876 & 2.833 & \azul{2.828} & 2.88 & \azul{2.88} & 0.875 & 0.89 & 0.386 & 0.4  \\    
$\et$ (deg.) & 0.055 & 0 & \azul{0.059} & 0.055 & \azul{0.053} & 0 & 0 &  0 & 0  \\    
\addlinespace[1pt]\cmidrule{1-10}
\end{tabular}
\label{tab.cone}
\end{table*}

\section{Conclusions}
\noindent
We have proposed in this work two fan-beam geometry (center of rotation) estimation algorithms and expanded them with a variable projection optimization approach to \azul{the problem of }geometry estimation in cone-beam tomography (horizontal detector shift and in-plane detector rotation). The two fan-beam techniques, based on a well known symmetry condition of fan-beam sinograms, are low-cost compared to reconstruction-based methods and have been proved to outperform state-of-the-art low-cost methods. The cone-beam geometry estimation is proved to be competitive to a method based on scanning a reference object previously manufactured for the two mentioned geometrical parameters. It is a gradient descent algorithm with Armijo backtracking line search for the step size at each iteration. Results are validated with simulated and industrial CT data with available python codes to reproduce the simulation experiments and to be adapted to any experimental fan- or cone-beam data.

\subsection*{Acknowledgements}
\noindent PG was supported by FWO Senior Research Project G0C2423N. SB was supported by the VLAIO Baekeland mandate grant HBC.2020.2280. RS was supported by the EU H2020-MSCA-ITN-2020 xCTing Project 956172. 

\bibliographystyle{abbrv} 
\bibliography{main}
\end{document}